\documentclass[12pt,leqno]{article}
\usepackage{amssymb}
\usepackage[mathscr]{eucal}
\usepackage{amsmath,amssymb,latexsym,theorem,bbm}
\usepackage{amsmath,amssymb,latexsym,theorem,enumerate}
\usepackage{color,url}
\usepackage{appendix}

\newcommand{\DD}{\mathbb{D}}

\newcommand{\NN}{\mathbb{N}}

\newcommand{\RR}{\mathbb{R}}

\newcommand{\ZZ}{\mathbb{Z}}

\newcommand{\bm}{{\boldsymbol{m}}}

\newcommand{\bq}{{\boldsymbol{q}}}
\newcommand{\bQ}{{\boldsymbol{Q}}}

\newcommand{\bx}{{\boldsymbol{x}}}

\newcommand{\by}{{\boldsymbol{y}}}

\newcommand{\bxi}{{\boldsymbol{\xi}}}

\newcommand{\bzero}{{\boldsymbol{0}}}
\newcommand{\bone}{{\boldsymbol{1}}}

\newcommand{\cA}{{\mathcal A}}
\newcommand{\cB}{{\mathcal B}}
\newcommand{\cC}{{\mathcal C}}
\newcommand{\cD}{{\mathcal D}}

\newcommand{\cL}{{\mathcal L}}
\newcommand{\cM}{{\mathcal M}}

\newcommand{\cN}{{\mathcal N}}

\newcommand{\cP}{{\mathcal P}}

\newcommand{\cW}{{\mathcal W}}

\newcommand{\dd}{\mathrm{d}}
\newcommand{\ee}{\mathrm{e}}

\newcommand{\EE}{\operatorname{\mathbb{E}}}

\newcommand{\PP}{\operatorname{\mathbb{P}}}

\newcommand{\cov}{\operatorname{Cov}}

\renewcommand{\leq}{\leqslant}
\renewcommand{\geq}{\geqslant}

\newcommand{\distr}{\stackrel{\cD}{\longrightarrow}}
\newcommand{\distre}{\stackrel{\cD}{=}}

\newcommand{\as}{\stackrel{{\mathrm{a.s.}}}{\longrightarrow}}

\newcommand{\bbone}{\mathbbm{1}}

\newcommand{\proofend}{\hfill\mbox{$\Box$}}

\numberwithin{equation}{section}

\theoremstyle{change} \theorembodyfont{\em}
\newtheorem{Lem}{Lemma.}[section]
\newtheorem{Thm}[Lem]{Theorem.}
\newtheorem{Pro}[Lem]{Proposition.}
\newtheorem{Cor}[Lem]{Corollary.}
\newtheorem{Def}[Lem]{Definition.}

\theorembodyfont{\rm}
\newtheorem{Rem}[Lem]{Remark.}
\newtheorem{Ex}[Lem]{Example.}

\begin{document}

\begin{center}
 {\bfseries\Large  Random means generated by random variables:\\[1mm] expectation and limit theorems}

\vspace*{3mm}

{\sc\large
  M\'aty\'as $\text{Barczy}^{*,\diamond}$,
  \ P\'al $\text{Burai}^{**}$ }

\end{center}

\vskip0.2cm

\noindent
 * MTA-SZTE Analysis and Stochastics Research Group,
   Bolyai Institute, University of Szeged,
   Aradi v\'ertan\'uk tere 1, H--6720 Szeged, Hungary.

\noindent
 ** Department of Differential Equations, Institute of Mathematics,
    Faculty of Natural Sciences, Budapest University of Technology and Economics,
    Egry J\'ozsef utca 1, H--1111 Budapest, Hungary.

\noindent e-mail: barczy@math.u-szeged.hu (M. Barczy),
                  buraip@math.bme.hu  (P. Burai).

\noindent $\diamond$ Corresponding author.

\vskip0.2cm


\renewcommand{\thefootnote}{}

\vspace*{-10mm}

\begin{abstract}
We introduce the notion of a random mean generated by a random variable and give a construction of its expected value.
We derive some sufficient conditions under which strong laws of large numbers and some limit theorems
  hold for random means generated by the elements of a sequence of independent and identically distributed random variables.
\end{abstract}

\section{Introduction}

The theory of means (also called averages) is an important, rich and growing field of mathematics,
 and it has several applications in practice as well.
For a recent monograph on averaging functions and their applications, see Beliakov et al. \cite{Belaikovetal2016}.
Random-valued mappings (functions) also appear in various fields of mathematics such
 as in transportation theory (see, e.g., Panaretos and Zemel \cite{PanZem}),
 in iterative functional equations (see, e.g., Baron and Jarczyk \cite{BarJar}, Baron \cite{Bar},
 and Jarczyk and Jarczyk \cite{Jarczyk2018}) or in theory of random  measures (see, e.g., Kallenberg \cite[Chapter 1]{Kal}).

In the present paper we introduce random (valued) means generated by a random variable and we give a construction of their expectations as well,
 which turn out to be usual (non-random) means.
Further, we derive some sufficient conditions under which strong laws of large numbers and some limit theorems
 hold for random means generated by the elements of a sequence of independent and identically distributed random variables.

 Concerning the notion of a random mean generated by a random variable given in Definition \ref{D:random_mean},
 there are (at least) two related notions in the literature, namely, the (continuous) family of means in the sense of
 P\'ales and Zakaria \cite[page 794]{PalZak} and the (continuous) random mean in the sense of Jarczyk and Jarczyk \cite[page 6838]{Jarczyk2020}.
In Remark \ref{Rem_JarJar_comparison} we compare our Definition \ref{D:random_mean} with these two related concepts.
Here we only note that in their definitions there is no random variable a priori involved, and in their setups
 a random mean is a real valued mapping, while our random mean generated by a random variable maps into the set of
 continuous (non-random) means on a nondegenerate, compact interval of the real numbers.
For historical fidelity, we mention that our Definition \ref{D:random_mean} was motivated by the definition of
 a random mean due to Jarczyk and Jarczyk \cite[page 6838]{Jarczyk2018}.

In Section \ref{Sec_random_mean} we introduce a notion of a random (valued) mean generated by a random variable.
Roughly speaking, given an $\RR^d$-valued random variable $\xi$, a $p$-variable random mean generated by $\xi$
 is a measurable function $M$ from a probability space to the space consisting of all (continuous) means on a nondegenerate, compact interval $I$
 of the real numbers $\RR$ such that there exists an auxiliary Borel measurable mapping $M_\xi:I^p\times\RR^d\to I$ in a way that the
 mappings $M_\xi(\cdot,\xi(\omega))$ and $(M(\omega))(\cdot)$ coincides for almost every $\omega$.
For a precise definition, see Definition \ref{D:random_mean}.
 We illustrate this definition by presenting a general method for constructing such random means (see Theorem \ref{Pro_3}),
 and we also give some examples such as  discrete random means generated by discrete random variables and
 random H\"older means (see Examples \ref{Ex_1} and \ref{Ex_2}).

The concept of Bochner integral (integral of maps defined on a measure space with values in a Banach space)
 allows us to define the expectation of a random mean generated by a random variable, see Definition \ref{Def_random_mean_exp}.
It turns out that the expectation in question is a usual (non-random) mean, see Theorem \ref{Pro_random_mean_exp}.
In Examples \ref{Ex_3} and \ref{Ex_4} we calculate the expectation of some random means generated by random variables given in Examples \ref{Ex_1} and \ref{Ex_2}.
We derive that the expectation of a $2$-variable random H\"older mean with weights governed by a uniform distribution on the interval $(0,1)$
 is nothing else but a Cauchy mean corresponding to some power functions or a logarithmic mean, see Example \ref{Ex_4}.
 In Remark \ref{Rem_Future_work}, motivated by Examples \ref{Ex_2} and \ref{Ex_4}, we initiate some possible future research directions.

Concerning the expectation of random means generated by a random variable, in probability theory there exist
 (at least) two somewhat related notions, namely, the expectation (also called barycenter)
 of a random probability measure on a compact metric space, see, e.g., Borsato et al.\ \cite[Appendix A.2]{BorHorRig}, and
 the Fr\'echet mean of a random measure with values in the \ $2$-Wasserstein space on \ $\RR^d$, \
 see, e.g., Panaretos and Zemel \cite[Section 3.2]{PanZem}.
In Remark \ref{Rem_random_mean_comparison} we recall both notions in order to see the similarities and differences
 with the expectation of a random mean generated by a random variable given in Definition \ref{Def_random_mean_exp}.

In Section \ref{Section_Limit} we derive some sufficient conditions under which strong laws of large numbers and some limit theorems
 hold for random means generated by the elements of a sequence of independent and identically distributed random variables.
More precisely, given a sequence of random means generated by the elements of a sequence of independent and identically distributed random variables,
 we consider the supremum norm of the difference of the arithmetic mean of the first $n$ random means in question and the common expectation of the random means,
 and we investigate the asymptotic behaviour of this random quantity as \ $n\to\infty$.
\ In Theorem \ref{Pro_SLLN_CLT}, the underlying sequence consists of independent and identically distributed discrete random variables,
 but the range of the random means in question is arbitrary in the sense that we do not suppose any special form of the usual (non-random) means
 in the range.
In Corollary \ref{Cro_Bernoulli} we consider a special case of Theorem \ref{Pro_SLLN_CLT}, namely,
 when the underlying sequence consists of  independent, identically and Bernoulli distributed random variables
 and the range of the random means in question is the set consisting of the arithmetic and geometric means in \ $[0,1]$.
In  Theorem \ref{Pro_uniform} we  establish some limit theorems for $p$-variable randomly weighted arithmetic means, when
 the underlying independent and identically distributed random variables are not necessarily discrete,
 so this result is out of scope of Theorem \ref{Pro_SLLN_CLT}.
In case of \ $p = 2$ \ we also formulate a corollary of Theorem \ref{Pro_SLLN_CLT} by simplifying the limit distribution,
 see Corollary \ref{Cor_uniform}.
Finally, we provide limit theorems for randomly weighted power means (which can be also called random H\"older means),
 and in this case instead of the arithmetic mean of the first $n$ random means in question we consider their geometric mean,
 and so the limit theorems have somewhat different forms compared to the previous ones.
 For a comparison of our limit theorems for random means with the Wasserstein law of large numbers for Fr\'echet means, see Remark \ref{Rem_Frechet}.

Section \ref{Section_proofs} is devoted to the proofs of Section \ref{Section_Limit}.
The main tools are the Kolmogorov's strong law of large numbers
 and the multidimensional central limit theorem together with the continuous mapping theorem.

We close the paper with an appendix, where we recall and prove a result on the continuity of the supremum for a two-variable continuous
 real-valued function by taking the supremum in one of its variables, see Theorem \ref{Thm_continuity}.
This result is used in the applications of the continuous mapping theorem in some of the proofs in Section \ref{Section_proofs}.

\section{Random means and their expectation}\label{Sec_random_mean}

Let \ $\NN$, \ $\ZZ_+$, \ $\RR$ \ and \ $\RR_+$ \ denote the sets of positive integers, non-negative integers, real numbers and non-negative real numbers,
respectively.
An interval $I\subset\RR$  is called nondegenerate if it contains at least two distinct points.
We denote by \ $\langle \bx,\by\rangle$ \ the Euclidean inner product of \ $\bx,\by\in\RR^d$, \  where \ $d\in\NN$.
 \ The Borel \ $\sigma$-algebra on \ $\RR^d$ \ is denoted by \ $\cB(\RR^d)$.
\ Convergence almost surely, convergence in distribution and equality in distribution will be denoted by \ $\as$, \ $\distr$ \ and \ $\distre$, \ respectively.
For any \ $d\in\NN$, \ $\cN_d(\bzero,\bQ)$ \ denotes a \ $d$-dimensional normal distribution on \ $\RR^d$ \ with mean vector \ $\bzero\in\RR^d$ \
 and covariance matrix \ $\bQ\in\RR^{d\times d}$.
\ In case of \ $d=1$, \ instead of \ $\cN_1$ \ we simply write \ $\cN$, \ and \ $\cN(0,0)$ \ denotes the Dirac distribution concentrated at $0$,
 \ i.e., for each \ $A\in\cB(\RR)$, \ $(\cN(0,0))(A):=1$ \ if \ $0\in A$, \ and \ $(\cN(0,0))(A):=0$ \ if \ $0\notin A$.

Let $I$ be a nondegenerate, compact interval of \ $\RR$, \ $p\in\NN$ \ be a positive integer,
 and let us denote by $\mathcal{C}(I^p)$ the vector space of real-valued,
 continuous functions defined on $I^p$, which becomes a Banach space with the usual supremum norm  given by
 \ $\|u\|:=\sup_{x\in I^p}|u(x)|$ \ for \ $u\in \mathcal{C}(I^p)$.
 The Borel \ $\sigma$-algebra on \ $\mathcal{C}(I^p)$ \ is denoted by \ $\cB(\mathcal{C}(I^p))$.
\  Given a function \ $f:I^p\times \RR^d \to I$, \ for any \ $y\in\RR^d$, \ we will denote by \ $f(\cdot,y)$ \ the function \ $I^p\ni x \mapsto f(x,y)$.

An $m\in\mathcal{C}(I^p)$ is said to be a $p$-variable, continuous mean on $I$ if
\[
\min(x)\leq m(x)\leq\max(x),\qquad x\in I^p,
\]
where \ $\min(x)$ \ and \ $\max(x)$ \ denotes the minimum and the maximum of the coordinates of \ $x\in I^p$, \ respectively.
If the above inequalities are strict whenever $x$ has at least two different coordinates, then $m$ is called  a strict mean.

From now on, we just simply use the terminology mean instead of continuous mean on $I$ if there is no ambiguity.
If \ $p=1$, \ then the only \ $1$-variable mean \ $m$ \ is \ $m(x)=x$, \ $x\in I$.

The subset $\cM_p\subset\mathcal{C}(I^p)$ denotes the class of $p$-variable means.

\begin{Pro}\label{P:Mp_convex}
The set \ $\cM_p\subset\mathcal{C}(I^p)$ \ is convex, bounded and closed with respect to the supremum norm.
\end{Pro}

\noindent{\bf Proof.}
The convexity of \ $\cM_p$ \ follows easily.
Concerning boundedness, if $m\in \cM_p$, then using that \ $m(x_1,\ldots,x_p) = x_1$ \ whenever \ $x_1 = \cdots = x_p \in I$, \ we have
 \[
 \| m\| =\sup_{x\in I^p}|m(x)| =  \max_{t\in I} \vert t \vert <\infty,
 \]
since $I$ is compact.
This means that the elements of $\cM_p$ have the same  (supremum) norm.

For proving closedness, let us  assume that $m_n\in \cM_p,\ n\in\NN$, is a convergent sequence in $\mathcal{C}(I^p)$,
 and let us denote by $u\in \mathcal{C}(I^p)$ its limit.
We intend to prove that
\begin{align}\label{help15}
\min(x)\leq u(x)\leq\max(x),\qquad x\in I^p.
\end{align}
Because convergence in norm implies pointwise convergence, we have
\begin{align}\label{help16}
m_n(x)\to u(x) \qquad \mbox{ as }n\to\infty, \qquad x\in I^p.
\end{align}
 Since \ $m_n\in\cM_p$, \ $n\in\NN$, \ we have \ $\min(x)\leq m_n(x)\leq\max(x)$, \ $x\in I^p$, \ $n\in\NN$.
\ So using \eqref{help16}, we have \eqref{help15}, as desired.
\proofend

For \ $m_1,m_2\in \cM_p$, \ let
 \begin{align}\label{Rem_metric}
    \varrho(m_1,m_2):=\sup_{x\in I^p} \vert m_1(x) - m_2(x) \vert.
 \end{align}
Note that $\cM_p$ is metric space furnished with $\varrho$ as a metric, but not a linear space.

\begin{Def}\label{D:random_mean}
Let $(\Omega,\cA,\PP)$ be a probability space and $\xi\colon\Omega\to\RR^d$ be a $d$-dimensional random variable, where $d\in\NN$. The map $M\colon\Omega\to\cM_p$ is called a ($p$-variable, continuous) random mean (in $I$) generated by $\xi$ if the following conditions are fulfilled:
\begin{enumerate}[(i)]
\item $M$ is $(\cA,\cB(\mathcal{C}(I^p))$-measurable, that is to say, \ $M$ \ is an \ $\cM_p$-valued random variable,
\item there is a \ $(\cB(I^p)\times\cB(\RR^d),\cB(I))$-measurable map \ $M_\xi\colon I^p\times\RR^d\to I$ \ such that%
 \begin{enumerate}[(a)]
  \item for every \ $\omega\in\Omega$ \ the map
        \ $I^p\ni x\mapsto M_\xi(x,\xi(\omega))$ \ is in \ $\mathcal{C}(I^p)$, \ and
   \item $\PP\left(\{\omega\in\Omega :(M(\omega))(x)=M_\xi(x,\xi(\omega))\ \mbox{ for every} \ x\in I^p\}\right)=1$.
\end{enumerate}
\end{enumerate}
\end{Def}

 Concerning the notion of a random mean generated by a random variable given in Definition \ref{D:random_mean},
 there are (at least) two related notions in the literature, namely, the (continuous) family of means in the sense of
 P\'ales and Zakaria \cite[page 794]{PalZak} and the (continuous) random mean in the sense of Jarczyk and Jarczyk \cite[page 6838]{Jarczyk2020}.
In the next remark we compare Definition \ref{D:random_mean} with these two related  concepts.

\begin{Rem}\label{Rem_JarJar_comparison}
\begin{enumerate}[(i)]
\item If \ $M$ \ is a random mean generated by \ $\xi$ \ in the sense of Definition \ref{D:random_mean} such that
 the map \ $M_\xi:I^p\times \RR^d \to I$ \ (appearing in part (ii) of Definition \ref{D:random_mean}) satisfies
 that for every \ $y\in\RR^d$ \ the map \ $I^p \ni x\mapsto M_\xi(x,y)$ \ is in \ $\cM_p$, \ then
 for each non-empty open interval \ $J$ \ of \ $I$, \ we have \ $M_\xi$ \ restricted to \ $J^p\times \RR^d$ \ is a continuous family
 of $p$-variable means on \ $J$ \ corresponding to the measurable space \ $(\RR^d, \cB(\RR^d))$ \ in the sense of
 P\'ales and Zakaria \cite[page 794]{PalZak}.
\item If \ $M$ \ is a random mean generated by \ $\xi$ \ in the sense of Definition \ref{D:random_mean}
 such that \ $M_\xi:I^p\times \RR^d \to I$ \ (appearing in part (ii) of Definition \ref{D:random_mean})
 satisfies that for every \ $y\in\RR^d$ \ the map \ $I^p \ni x\mapsto M_\xi(x,y)$ \ is in \ $\cM_p$,
 then $M_\xi$ is a continuous random mean on  $I$  in the sense of Jarczyk and Jarczyk \cite[page 6838]{Jarczyk2020}
 corresponding to the probability space $(\RR^d,\cB(\RR^d), \PP_\xi)$,
 where $\PP_\xi$ denotes the distribution of $\xi$ (i.e., \ $\PP_\xi(A) := \PP(\xi\in A)$, \ $A\in\cB(\RR^d)$).
Indeed,  $M_\xi$ is $(\cB(I^p) \times \cB(\RR^d), \cB(I))$-measurable, \ $\{ y\in\RR^d:  M_\xi( \cdot , y)\in \cM_p \}=\RR^d$ \ and hence
 \[
 \PP_\xi( \{ y\in\RR^d:  M_\xi( \cdot , y)\in \cM_p \}) = \PP_\xi(\RR^d) = 1.
 \]
\end{enumerate}

We point out to the facts that both in the definition of a (continuous) family of $p$-variable means on $J$ due to
 P\'ales and Zakaria \cite[page 794]{PalZak} and in the definition of a (continuous) random mean  on $I$ due to Jarczyk and Jarczyk \cite[page 6838]{Jarczyk2020}
 there is no random variable a priori involved.

 Further, a (continuous) family of means on $J$ and a random mean on $I$ is a $J$-valued
 and an $I$-valued mapping, respectively, while our random mean generated by a random variable is an \ $\cM_p$-valued mapping.
\proofend
\end{Rem}

Next, we illustrate the definition of a random mean generated by a random variable by presenting a general method for constructing
 such random means, and we also give some examples.

\begin{Thm}\label{Pro_3}
Let $(\Omega,\cA,\PP)$ be a probability space, \ $\xi:\Omega\to\RR^d$ \ be a $d$-dimensional random variable
 and \ $f \colon I^p\times\RR^d\to I$ \ be a \ $(\cB(I^p)\times\cB(\RR^d),\cB(I))$-measurable map
 such that for each \ $y\in\RR^d$, \ the map \ $I^p\ni x\mapsto f(x,y)$ \ is in \ $\cM_p$.
Then the map \ $M:\Omega\to \cM_p$ \ given by
 \[
   (M(\omega))(x) := f(x,\xi(\omega)), \qquad \omega\in\Omega, \;\; x\in I^p,
 \]
 is a random mean generated by \ $\xi$.
 \end{Thm}

\noindent{\bf Proof.}
By the construction, \ $M(\omega)\in\cM_p$ \ for all \ $\omega\in\Omega$, \ and part (ii) of Definition \ref{D:random_mean} holds
 with the choice $M_\xi:=f$.

Further, \ $M$ \ can be written as \ $M=\varphi \circ \xi$, \ where \ $\varphi : \RR^d\to\cM_p$, \ $\varphi(y):=f(\cdot,y)$, \ $y\in\RR^d$.
\ Here \ $\xi$ \ is \ $(\cA,\cB(\RR^d))$-measurable and  we check that \ $\varphi$ \ is \ $(\cB(\RR^d), \cB(\cC(I^p))$-measurable.
It is known that \ $\cB(\cC(I^p))$ coincides with the $\sigma$-algebra generated by $\cC$,
 where \ $\cC$ \ is the set of so-called cylinder sets of \ $\cC(I^p)$ \ having the form
 \[
   \big\{ g\in \cC(I^p) : (g(t_1), ..., g(t_n))\in B_1\times \cdots\times B_n\big\},
   \qquad n\in\NN, \;\; t_1,\ldots,t_n\in I^p, \;\; B_1,\ldots,B_n\in\cB(\RR),
 \]
 see, e.g., Kuo \cite[Chapter I, Theorem 4.2]{Kuo1975}.
So it is enough to check that for  all \ $n\in\NN$, \ $t_1,\ldots,t_n\in I^p$ \ and \ $B_1,\ldots,B_n\in\cB(\RR)$, \ we have
 \[
   \varphi^{-1}\Big( \big\{ g\in \cC(I^p) : (g(t_1), ..., g(t_n))\in B_1\times \cdots\times B_n\big\} \Big)\in\cB(\RR^d).
 \]
Here
 \begin{align*}
   &\varphi^{-1}\Big( \big\{ g\in \cC(I^p) : (g(t_1), ..., g(t_n))\in B_1\times \cdots\times B_n\big\} \Big)\\
   & = \Big\{ y\in\RR^d : f(\cdot, y)\in\cC(I^p) \;\;\text{and} \;\; (f(t_1,y),\ldots,f(t_n,y))\in B_1\times\cdots\times B_n \Big\}\\
   & = \Big\{ y\in\RR^d : f(\cdot, y)\in\cC(I^p) \Big\}\cap \bigcap_{i=1}^n f_{t_i}^{-1}(B_i),
 \end{align*}
 where for each \ $i=1,\ldots,n$, \ the function \ $f_{t_i}: \RR^d\to I$, \ $f_{t_i}(y):=f(t_i,y)$, \ $y\in\RR^d$, \ is a section of \ $f$ \ on \ $\RR^d$.
\ By our assumptions,
\[
 \big\{ y\in\RR^d : f(\cdot, y)\in\cC(I^p) \big\}
  = \big\{ y\in\RR^d : f(\cdot, y)\in\cM_p \big\}
   = \RR^d,
 \]
 and \ $f_{t_i}$, \ $i=1,\ldots,n$, \ is \ $(\cB(\RR^d), \cB(I))$-measurable
 \ (see, e.g., Cohn \cite[Lemma 5.1.2]{Cohn1980}) yielding that \ $f_{t_i}^{-1}(B_i)\in\cB(\RR^d)$, \ $i=1,\ldots,n$.
\ As a consequence, \ $\varphi^{-1}(C)\in \cB(\RR^d)$ \ for all \ $C\in\cC$, \ as desired.
Consequently, \ $M$ \ is \ $(\cA,\cB(\cC(I^p)))$-measurable yielding part (i) of Definition \ref{D:random_mean}.
\proofend

\begin{Ex}[Discrete random mean generated by a discrete random variable]\label{Ex_1}
Let \ $\xi$ \ be a  (one-dimensional) discrete random variable.
In this example without loss of generality, we can and do assume that the range of \ $\xi$ \ is \ $\ZZ_+$.
\ Let us consider a sequence \ $m_i$, \ $i\in\ZZ_+$, \ in \ $\cM_p$,
 \ and let \ $M:\Omega\to \cM_p$,
 \[
     (M(\omega))(x) := \sum_{i=0}^\infty m_i(x) \bbone_{\{i\}}(\xi(\omega))   = m_{\xi(\omega)}(x),
       \qquad \omega\in\Omega, \;\; x\in I^p.
 \]
We check that \ $M$ \ is a random mean generated by \ $\xi$.
\ First, note that \ $M(\omega)\in\cM_p$ \ for all \ $\omega\in \Omega$, \ since  $\xi$ is $\ZZ_+$-valued and
 if \ $\omega\in\Omega$ \ is such that \ $\xi(\omega) = i$, \ where \ $i\in\ZZ_+$, \ then \ $(M(\omega))(x) = m_i(x)$, \  $x\in I^p$, \ and $m_i\in\cM_p$.

Next, we check that \ $M$ \ is \ $(\cA,\cB(\cC(I^p)))$-measurable.
For each \ $n\in\NN$, \ let
 \[
 \left( M^{(n)}(\omega) \right)(x):=\sum_{i=0}^n m_i(x)  \bbone_{\{i\}}(\xi(\omega)) 
                   = \sum_{i=0}^n m_i(x) \bbone_{\xi^{-1}(\{i\})}(\omega), \qquad \omega\in\Omega.
 \]
Then \ $M^{(n)}$ \ is a simple \ $\cM_p$-valued random variable  (i.e., it has only finitely many values)
 being \ $(\cA,\cB(\cC(I^p)))$-measurable, and for each
 \ $\omega\in\Omega$, \ we have
 \begin{align*}
  &\sup_{x\in I^p} \vert (M^{(n)}(\omega))(x) - (M(\omega))(x) \vert
   = \sup_{x\in I^p} \left\vert \sum_{i=n+1}^\infty m_i(x)  \bbone_{\{i\}}(\xi(\omega))  \right\vert\\
  &\qquad  = \sup_{x\in I^p} \left\vert m_{\xi(\omega)}(x)  \bbone_{[n+1,\infty)}(\xi(\omega))   \right\vert
           = \bbone_{[n+1,\infty)}(\xi(\omega))  \sup_{x\in I^p} \vert m_{\xi(\omega)}(x)\vert
    \to 0 \qquad \text{as \ $n\to\infty$,}
 \end{align*}
 since \ $\sup_{x\in I^p} \vert m_{\xi(\omega)}(x)\vert<\infty$ \  (as we have seen in the proof of Proposition \ref{P:Mp_convex}).
\ This yields that \ $M$ \ is a pointwise limit of \ $M^{(n)}$ \ as \ $n\to\infty$, \ so \ $M$ \ is  \ $(\cA,\cB(\cC(I^p)))$-measurable as well.
So part (i) of Definition \ref{D:random_mean} holds.

Let \ $M_\xi:I^p\times \RR\to I$ \ be given by
 \[
    M_\xi(x,y):= \sum_{i=0}^\infty m_i(x) \bbone_{\{i\}}(y),
    \qquad x\in I^p,\;\; y\in\RR.
 \]
For each \ $n\in\NN$, \ let \ $M_\xi^{(n)}:I^p\times \RR\to I$ \ be given by
 \[
    M_\xi^{(n)}(x,y):= \sum_{i=0}^n m_i(x) \bbone_{\{i\}}(y),
    \qquad x\in I^p,\;\; y\in\RR.
 \]
Then \ $M_\xi^{(n)}$ \ is \ $(\cB(I^p)\times \cB(\RR), \cB(I))$-measurable, since \ $ m_i$ \ is continuous and $\{i\}\in\cB(\RR)$
 for each $i\in\ZZ_+$.
\ Further, for all \ $x\in I^p$ \ and \ $y\in\RR\setminus \ZZ_+$, \ we have \ $\vert M_\xi^{(n)}(x,y) - M_\xi(x,y) \vert = \vert 0-0\vert=0$, \ and
 for all \ $x\in I^p$ \ and \ $y\in\ZZ_+$, \ we have
 \begin{align*}
  \vert M_\xi^{(n)}(x,y) - M_\xi(x,y) \vert
    = \left\vert\sum_{i=n+1}^\infty m_i(x) \bbone_{\{i\}}(y) \right\vert
    = \vert m_y(x)\vert \bbone_{[n+1,\infty)}(y)
     \to 0\qquad \text{as \ $n\to\infty$.}
 \end{align*}
So \ $M_\xi$ \ is a pointwise limit of \ $M_\xi^{(n)}$ \ as \ $n\to\infty$, \ yielding that \ $M_\xi$ \ is
 \ $(\cB(I^p)\times \cB(\RR), \cB(I))$-measurable.
Moreover, for all \ $\omega\in\Omega$,
 \[
  M_\xi(x,\xi(\omega)) = \sum_{i=0}^\infty m_i(x)  \bbone_{\{i\}}(\xi(\omega)) 
                       = (M(\omega))(x),\qquad x\in I^p.
 \]
At the beginning of the example we have seen that $M(\omega)\in\cM_p$ for all $\omega\in\Omega$.
So, we get that part (ii) of Definition \ref{D:random_mean} holds as well.

 We can call \ $M$ \ a discrete random mean generated by the discrete random variable \ $\xi$ \ in question,
 since the range of \ $M$ \ contains countably many elements of $\cM_p$.
\proofend
\end{Ex}

\begin{Ex}[Random H\"older means]\label{Ex_2}
 If \ $I:=[a,b]$, \ where \ $0<a<b<\infty$, \ then let \ $f:I^2\times \RR \times(0,1)\to I$  \ be defined by
 \[
  f(x_1,x_2,\alpha,\lambda)
   :=\begin{cases}
       \left(\lambda x_1^{\alpha}+(1-\lambda)x_2^{\alpha}\right)^{\frac{1}{\alpha}},& \mbox{ if }\ \alpha\not=0,\\
 x_1^{\lambda}x_2^{1-\lambda},&\mbox{ if }\ \alpha=0,
 \end{cases}
 \qquad (x_1,x_2,\alpha,\lambda)\in I^2\times\RR\times(0,1),
 \]
 and let \ $\xi\colon\Omega\to\RR\times(0,1)$ \ be a random variable.
If \ $I:=[0,b]$, \ where \ $0<b<\infty$, \ then let \ $f:I^2\times \RR_+ \times(0,1)\to I$  \ be defined by
 \[
  f(x_1,x_2,\alpha,\lambda)
   :=\begin{cases}
       \left(\lambda x_1^{\alpha}+(1-\lambda)x_2^{\alpha}\right)^{\frac{1}{\alpha}},& \mbox{ if }\ \alpha>0,\\
 x_1^{\lambda}x_2^{1-\lambda},&\mbox{ if }\ \alpha=0,
 \end{cases}
 \qquad (x_1,x_2,\alpha,\lambda)\in I^2\times\RR_+\times(0,1),
 \]
 and let \ $\xi\colon\Omega\to\RR_+\times(0,1)$ \ be a random variable.
Note that for each \ $(\alpha,\lambda)\in\RR\times(0,1)$ \ in case of \ $I=[a,b]$ \ ($0<a<b<\infty$),
 and for each \ $(\alpha,\lambda)\in\RR_+\times(0,1)$ \ in case of \ $I=[0,b]$ \ ($0<b<\infty$),
 \ the map \ $I^2\ni (x_1,x_2)\mapsto f(x_1,x_2,\alpha,\lambda)$ \ is a H\"older mean (also called weighted power mean), so it is in \ $\cM_2$. \
Hence, using also that \ $f$ \ is \ $(\cB(I^2)\times\cB(\RR^2),\cB(I))$-measurable,
 we can apply Theorem \ref{Pro_3} and we have \ $M:\Omega\to\cM_2$, \ $(M(\omega))(x_1,x_2):=f(x_1,x_2,\xi(\omega))$, \ $\omega\in\Omega$,
 \ $(x_1,x_2)\in I^2$, \ is a random mean generated by \ $\xi$, \ which can be called a random H\"older mean.
One can define the $p$-variable version of this random mean in a similar way.
\proofend
\end{Ex}

The next proposition allows us to define the expected value of a random mean generated by a random variable $\xi$.
The concept of Bochner integrability, integral of maps defined on a measure space with values in a Banach space, has a key role. We use the results and terminology of  Cohn \cite[Appendix E]{Cohn1980}.

\begin{Pro}\label{proposition:bochner_integrability}
If $M\colon\Omega\to \cM_p$ is $(\cA,\cB(\mathcal{C}(I^p)))$-measurable, then it is Bochner integrable.
\end{Pro}

\noindent{\bf Proof.}
The function $\Omega\ni\omega\mapsto M(\omega)$  (considered as a function with values in \ $\cC(I^p)$) \
 is Bochner integrable if it  is strongly measurable -- i.e., $M$ is $(\cA,\cB(\cC(I^p))$-measurable and has a separable range,
 where, by the range of $M$ we mean the subset $M(\Omega)$ of $\cC(I^p)$ -- and the function $\Omega\ni\omega\mapsto \|M(\omega)\|$ is integrable
 with respect to \ $\PP$, \ see Cohn \cite[Appendix E]{Cohn1980}.

Because of the Stone-Weierstrass approximation theorem, $\cC(I^p)$ is separable,
 and since each subspace of a separable metric space is separable, we have the range of $M$ is also separable.
Hence the $(\cA,\cB(\cC(I^p)))$-measurability of \ $M$ \ implies that \ $M$ \ is strongly measurable.

Moreover, using that \ $(M(\omega))(x_1,\ldots,x_p) = x_1$ \ whenever \ $x_1 = \cdots = x_p \in I$ \ and \ $\omega\in\Omega$ \
 (due to \ $M(\omega)\in\cM_p$, \ $\omega\in \Omega$), \ we have
 \begin{align}\label{help14}
  \|M(\omega)\|=\sup_{x\in I^p}| (M(\omega))(x)|=\max_{t\in I}|t|  <\infty , \qquad \omega\in\Omega,
 \end{align}
 since \ $I$ \ is compact.
 So the function \ $\Omega\ni\omega\mapsto \|M(\omega)\|$ \ is the constant  \ $\max_{t\in I}|t|$ \ function, and, using the fact that \ $\PP(\Omega)=1$,
 \ we have that it is integrable.
\proofend

\begin{Cor}\label{Cor:integrability_of_random_means}
If \ $M:\Omega\to \cM_p$ \ is a random mean generated by a $d$-dimensional random variable \ $\xi$, \ then it is Bochner integrable.
\end{Cor}

\noindent\textbf{Proof.}
Since \ $M$ \ is \ $(\cA,\cB(\cC(I^p)))$-measurable (following from part (i) of Definition \ref{D:random_mean}),
 Proposition \ref{proposition:bochner_integrability} yields the statement.
\proofend

According to Corollary \ref{Cor:integrability_of_random_means} the following definition does make sense.

\begin{Def}\label{Def_random_mean_exp}
Let $M:\Omega\to \cM_p$ be a random mean generated by a $d$-dimensional random variable $\xi$
  defined on a probability space $(\Omega,\cA,\PP)$.
Then the element of \ $\mathcal{C}(I^p)$ \ given by
 \[
 \EE(M):=\int_\Omega M(\omega)\,\dd\PP(\omega)
 \]
 is called the expected value or the expectation of $M$.
\end{Def}

\begin{Thm}\label{Pro_random_mean_exp}
If $M:\Omega\to \cM_p$ is a random mean generated by a $d$-dimensional random variable $\xi$,
 then its expected value \ $\EE(M)$ \ is a non-random mean, that is to say, $\EE(M)\in \cM_p$.
 \  Further,
  \[
   (\EE(M))(x) = \int_\Omega (M(\omega))(x)\,\dd\PP(\omega),
               \qquad x\in I^p.
  \]
\end{Thm}

\noindent{\bf Proof.}
By Corollary \ref{Cor:integrability_of_random_means}, \ $\EE(M)$ \ exists, and, especially, \ $\EE(M)\in \mathcal{C}(I^p)$.
So, it remains to check that it is in $\cM_p$.

It follows from Hyt\"onen et al. \cite[Proposition 1.2.12]{HytNerVerWei2016} that
 \[
   \EE(M)\in \overline{\mathrm{conv} \{M(\omega) : \omega\in\Omega\}},
 \]
 where \ $\mathrm{conv} \{M(\omega) : \omega\in\Omega\}$ \ denotes the convex hull of \ $\{M(\omega) : \omega\in\Omega\}$, \ and
  \ $\overline{\mathrm{conv} \{M(\omega) : \omega\in\Omega\}}$ \ is its closure in \ $\cC(I^p)$.
Additionally, Proposition \ref{P:Mp_convex} implies
\[
\overline{\mathrm{conv} \{M(\omega) : \omega\in\Omega\}}\subset \cM_p.
\]
 These two gives that $\EE(M)\in\cM_p$.

Further, for each \ $x\in I^p$, \ let \ $\varphi_x:\cC(I^p)\to\RR$, \ $\varphi_x(h):=h(x)$, \ $h\in\cC(I^p)$.
Then for each \ $x\in I^p$, \ $\varphi_x$ \ is a linear functional on \ $\cC(I^p)$, \ and hence
 Proposition E.11 in Cohn \cite{Cohn1980} yields that
 \begin{align}\label{help17}
  (\EE(M))(x) = \varphi_x(\EE(M)) =  \int_\Omega \varphi_x(M(\omega))\,\dd\PP(\omega)
              = \int_\Omega (M(\omega))(x)\,\dd\PP(\omega),
 \end{align}
as desired.

Using \eqref{help17} one can give another (a more elementary) proof of the fact that \ $\EE(M)\in\cM_p$.
\ Namely, since \ $M(\omega)\in\cM_p$, $\omega\in\Omega$, \ we have \ $\min(x)\leq (M(\omega))(x)\leq \max(x)$, $x\in I^p$, \
 thus
 \[
    \min(x) \leq \int_\Omega (M(\omega))(x)\,\dd\PP(\omega) \leq \max(x),\qquad x\in I^p,
 \]
 and, by \eqref{help17}, we have \ $\min(x) \leq (\EE(M))(x) \leq \max(x)$, \ $x\in I^p$, \ i.e., \ $\EE(M)\in \cM_p$.
\proofend

\begin{Rem}\label{Rem_random_mean_calculation}
Let $(\Omega,\cA,\PP)$ be a probability space, \ $\xi:\Omega\to\RR^d$ \ be a $d$-dimensional random variable
 and \ $f \colon I^p\times\RR^d\to I$ \ be a \ $(\cB(I^p)\times\cB(\RR^d),\cB(I))$-measurable map
 such that for each \ $y\in\RR^d$, \ the map \ $I^p\ni x\mapsto f(x,y)$ \ is in \ $\cM_p$.
Then, by Theorem \ref{Pro_3}, \ $M:\Omega\to \cM_p$ \ given by \ $(M(\omega))(x) := f(x,\xi(\omega))$,
 \ $\omega\in\Omega$, $x\in I^p$, is a random mean generated by \ $\xi$.
Further, by Theorem \ref{Pro_random_mean_exp},
 \begin{align}\label{help_random_mean_calculation}
 (\EE(M))(x) = \int_\Omega f(x,\xi(\omega)) \,\dd\PP(\omega)
               = \int_{\RR^d} f(x, y)\, \dd\PP_\xi(y),
               \qquad x\in I^p,
 \end{align}
 since the map $\RR^d\ni y \mapsto f(x,y)\in I$ is $(\cB(\RR^d),\cB(I))$-measurable for each fixed $x\in I^p$ and hence one can apply
 a result on integration with respect to an image measure (see, e.g., Cohn \cite[Proposition 2.6.8]{Cohn1980}).
Note also that in this case $\EE(M)$ depends only on $f$ and the distribution $\PP_\xi$ of $\xi$.
We do not know whether all the random means can be written in the form given in Theorem \ref{Pro_3}.
\proofend
\end{Rem}

 Next we recall the notions of expectation (also called barycenter) of a random probability measure on a compact metric space
 (see, e.g., Borsato et al.\ \cite[Appendix A.2]{BorHorRig}), and the Fr\'echet mean of a random measure with values in the \ $2$-Wassertein space
 on \ $\RR^d$ \ (see, e.g., Panaretos and Zemel \cite[Section 3.2]{PanZem}) in order to see the similarities and differences
 compared to the expected value of a random mean generated by a random variable given in Definition \ref{Def_random_mean_exp}.

\begin{Rem}\label{Rem_random_mean_comparison}
First, we recall the expectation of a random probability measure on a compact metric space.
Given a probability space \ $(\Omega,\cA,\PP)$ \ and a compact metric space \ $S$ \ endowed with its Borel
 \ $\sigma$-algebra \ $\cB(S)$, \ a random probability measure on \ $S$ \ is defined to be a Borel measurable map
 \ $\eta:\Omega\to \cP_1(S)$, \ where \ $\cP_1(S)$ \ denotes the set of probability measures on \ $(S,\cB(S))$ \ and
 \ $\cP_1(S)$ \ is endowed with the Borel $\sigma$-algebra corresponding to the topology of weak convergence
 according to which a sequence \ $(\mu_n)_{n\geq 1}$ \ in \ $\cP_1(S)$ \ converges to a given \ $\mu\in\cP_1(S)$ \
  if \ $\int_S f(s)\,\mu_n(\dd s) \to \int_S f(s)\,\mu(\dd s)$ \ as \ $n\to\infty$ \ for each continuous (hence bounded) function
  \ $f:S\to\RR$.
\ Then, as a consequence of Riesz-Markov's representation theorem,
 there exists a unique element \ $\EE(\eta)$ \ of \ $\cP_1(S)$ \ such that the equality
 \ $\int_S f(s) \, (\EE(\eta))(\dd s) = \int_\Omega \int_S f(s) \,\eta^{\omega}(\dd s)\PP(\dd \omega)$ \ holds for each
 continuous (hence bounded) function \ $f:S\to \RR$, \ where \ $\eta^\omega$ \ denotes the value of the random measure \ $\eta$ \
 at the point \ $\omega\in\Omega$, \ see, e.g., Borsato et al.\ \cite[Theorem A.6 and Definition A.7]{BorHorRig}.

Next, we recall the Fr\'echet mean of a random measure with values in the \ $2$-Wassertein space on \ $\RR^d$.
\ Given a probability space \ $(\Omega,\cA,\PP)$ \ and \ $d\in\NN$, \ the \ $2$-Wasserstein space on \ $\RR^d$ \ is defined by
 \[
     \cW_2(\RR^d) := \Big\{ \mu\in\cP_1(\RR^d) : \int_{\RR^d} \Vert \bx\Vert^2 \mu(\dd\bx)<\infty \Big\},
 \]
 where \ $\cP_1(\RR^d)$ \ denotes the set of probability measures on \ $\RR^d$.
\ For \ $\mu,\nu\in\cP_1(\RR^d)$, \ let \ $\Pi(\mu,\nu)$ \ be the set of probability measures
 \ $\pi\in\cP_1(\RR^d\times\RR^d)$ \ such that \ $\pi(A\times \RR^d) = \mu(A)$, \ $A\in\cB(\RR^d)$, \ and
 \ $\pi(\RR^d\times B) = \nu(B)$, \ $B\in\cB(\RR^d)$, \ i.e., \ $\mu$ \ and \ $\nu$ \ are the marginals of \ $\pi$.
\ The \ $2$-Wasserstein distance between \ $\mu$ \ and \ $\nu$ \ is defined as
 \[
     W_2(\mu,\nu):= \left(   \inf_{\pi\in\Pi(\mu,\nu)} \int_{\RR^d\times\RR^d} \Vert \bx_1 - \bx_2\Vert^2 \,\dd \pi(\bx_1,\bx_2)\right)^{\frac{1}{2}},
     \qquad \mu,\nu\in\cP_1(\RR^d).
 \]
Then \ $W_2$ \ is a metric on \ $\cW_2(\RR^d)$, \ see Villani \cite[Chapter 7]{Vil}.
By a random measure with values in \ $\cW_2(\RR^d)$, \ we mean a measurable map \ $\Lambda: \Omega\to\cW_2(\RR^d)$, \ where \ $\cW_2(\RR^d)$ \
 is endowed with its Borel $\sigma$-algebra (corresponding to the metric \ $W_2$).
\ By the Fr\'echet mean (expectation) of a random measure \ $\Lambda$ \ with values in \ $\cW_2(\RR^d)$, \ we mean the minimizer (if it is unique)
 of the Fr\'echet functional
 \[
    F(\gamma):=\frac{1}{2} \EE((W_2(\gamma,\Lambda))^2), \qquad \gamma\in \cW_2(\RR^d),
 \]
 see, e.g., Definition 3.2.1 in Panaretos and Zemel \cite{PanZem}.
We note that the Fr\'echet functional associated with any random measure \ $\Lambda$ \ with values in \ $\cW_2(\RR^d)$ \ admits a minimizer (see, e.g.,
 Panaretos and Zemel \cite[Proposition 3.2.3]{PanZem}), and for a result on the uniqueness of Fr\'echet means, see, e.g.,
 Proposition 3.2.7 in Panaretos and Zemel \cite{PanZem}.
 Further, see Remark \ref{Rem_Frechet} for a comparison of our forthcoming limit theorems for random means
 generated by random variables with the Wasserstein law of large numbers for Fr\'echet means (Panaretos and Zemel \cite[Corollary 3.2.10]{PanZem}).
\proofend
\end{Rem}

Next, we determine the expectation of the random means given in Examples \ref{Ex_1} and \ref{Ex_2}
 (in case of Example \ref{Ex_2} with special choices of \ $\xi$).

\begin{Ex}\label{Ex_3}
The expectation of the random mean generated by a discrete random variable \ $\xi$ \ having range in \ $\ZZ_+$
 \ given in Example \ref{Ex_1} takes the form
 \begin{align*}
  \EE(M)  = \int_\Omega M(\omega)\,\dd \PP(\omega)
           = \lim_{n\to\infty} \int_\Omega M^{(n)}(\omega)\,\dd \PP(\omega)
           = \lim_{n\to\infty} \sum_{i=0}^n m_i \PP(\xi=i)
           = \sum_{i=0}^\infty  \PP(\xi=i) m_i ,
 \end{align*}
 where  $M^{(n)}$ is introduced in Example \ref{Ex_1}, the series above converges in $\cC(I^p)$,
 and for the second equality we used the construction of Bochner integral (see, e.g., Cohn \cite[Appendix E]{Cohn1980}).
 \proofend
\end{Ex}

\begin{Ex}[Expectation of some random H\"older means]\label{Ex_4}
Let us consider the random H\"older mean \ $M$ \ given in Example \ref{Ex_2} generated by a random variable \ $\xi$.
\  First, let us suppose that the distribution \ $\PP_\xi$ \ of \ $\xi$ \ takes the form
 \ $\PP_\xi = \delta_{\alpha_0}\otimes \PP_U$, \ where \ $\alpha_0\in \RR_+$, \ $\delta_{\alpha_0}$
 \ denotes the Dirac measure concentrated at \ $\alpha_0$, \ and \ $U$ \ is a uniformly distributed random variable in the interval \ $(0,1)$.
 In Example \ref{Ex_2}, let us choose \ $I:=[0,b]$, \ where \ $b>0$. \\
 In case of \ $\alpha_0\in(0,\infty)$, for the expectation \ $\EE(M)\in\cM_2$ \ of \ $M$, \ we have
 \begin{align*}
   (\EE(M))(x_1,x_2) & = \int_{\RR^2} f(x_1,x_2,\alpha,\lambda)\,\PP_\xi(\dd\alpha,\dd\lambda)
                       = \int_0^1 (\lambda x_1^{\alpha_0} + (1-\lambda)x_2^{\alpha_0})^{\frac{1}{\alpha_0}}\,\dd\lambda\\
                     & = \int_0^1 \big(  (x_1^{\alpha_0} - x_2^{\alpha_0})\lambda+ x_2^{\alpha_0} \big)^{\frac{1}{\alpha_0}}\,\dd\lambda\\
     &=\begin{cases}
                        \frac{x_1^{\alpha_0+1} - x_2^{\alpha_0+1}}{(\frac{1}{\alpha_0} + 1)(x_1^{\alpha_0} - x_2^{\alpha_0})}&
                           \quad \text{if \ $x_1\ne x_2$,}\\
                        x_1 &\quad  \text{if \ $x_1=x_2$,}
     \end{cases}
                      \qquad x_1,x_2\in I,
 \end{align*}
  where  the first equality follows by \eqref{help_random_mean_calculation}.
\ In this case one can check that \ $\EE(M)$ \ is nothing else but a Cauchy mean corresponding to the power functions
 \ $x^{\alpha_0+1}$, \ $x\in I$, \ and \ $x^{\alpha_0}$, \ $x\in I$, \ see, e.g., Beliakov et al.\ \cite[Definition 2.50]{Belaikovetal2016}
  or Jarczyk and Jarczyk \cite[Section 5.1]{Jarczyk2018}.

In case of \ $\alpha_0 =0$, \ for the expectation \ $\EE(M)\in\cM_2$ \ of \ $M$, \ by \eqref{help_random_mean_calculation}, we have
 \begin{align*}
   (\EE(M))(x_1,x_2)  = \int_0^1 x_1^\lambda x_2^{1-\lambda}\,\dd\lambda
                      =\begin{cases}
                        \frac{x_1 - x_2}{\ln(x_1) - \ln(x_2)}&
                           \quad \text{if \ $x_1\ne x_2$, \ $x_1,x_2\in I\cap(0,\infty)$,}\\
                            x_1 &\quad  \text{if \ $x_1=x_2\in I\cap(0,\infty)$,}\\
                            0 &\quad  \text{if \ $x_1=0$ \ or \ $x_2=0$.}
                      \end{cases}
 \end{align*}
In this case, \ $\EE(M)$ \ restricted to \ $(0,\infty)\times(0,\infty)$ \  is nothing else but the logarithmic mean,
 see, e.g., Beliakov et al.\ \cite[Definition 2.45]{Belaikovetal2016}  or Jarczyk and Jarczyk \cite[Section 5.1]{Jarczyk2018}.

Next, let us suppose that the distribution \ $\PP_\xi$ \ of \ $\xi$ \ takes the form
 \ $\PP_\xi = \delta_{\alpha_0}\otimes \PP_V$, \ where \ $\alpha_0\in\RR_+$ \ and \ $V$ \ is a random
 variable with density function \ $f_V(\lambda) = 2\lambda \bone_{(0,1)}(\lambda)$, \ $\lambda\in\RR$.
\ As before, in Example \ref{Ex_2}, let us choose \ $I:=[0,b]$, \ where \ $b>0$. \\
In case of \ $\alpha_0\in(0,\infty)$, \ for the expectation \ $\EE(M)\in\cM_2$ \ of \ $M$, \ by \eqref{help_random_mean_calculation}, we have
 \begin{align*}
   &(\EE(M))(x_1,x_2)
      = \int_0^1 \big(  (x_1^{\alpha_0} - x_2^{\alpha_0})\lambda+ x_2^{\alpha_0} \big)^{\frac{1}{\alpha_0}} 2\lambda \,\dd\lambda\\
   & = \frac{2}{x_1^{\alpha_0} - x_2^{\alpha_0}} \int_0^1 \big(  (x_1^{\alpha_0} - x_2^{\alpha_0})\lambda+ x_2^{\alpha_0} \big)^{\frac{1}{\alpha_0} + 1} \,\dd\lambda
          - \frac{2x_2^{\alpha_0}}{x_1^{\alpha_0} - x_2^{\alpha_0}}
             \int_0^1 \big(  (x_1^{\alpha_0} - x_2^{\alpha_0})\lambda+ x_2^{\alpha_0} \big)^{\frac{1}{\alpha_0}} \,\dd\lambda\\
   & = \frac{2}{\frac{1}{\alpha_0} + 2} \frac{(x_1^{2\alpha_0+1} - x_2^{2\alpha_0+1})}{ (x_1^{\alpha_0} - x_2^{\alpha_0})^2 }
         - \frac{2}{\frac{1}{\alpha_0} + 1} \frac{x_2^{\alpha_0}(x_1^{\alpha_0+1} - x_2^{\alpha_0+1})}{ (x_1^{\alpha_0} - x_2^{\alpha_0})^2 } \\
  & = \frac{2}{\left( \frac{1}{\alpha_0} + 1\right) \left( \frac{1}{\alpha_0} + 2\right)}
        \cdot \frac{ \left( \frac{1}{\alpha_0} + 1\right) x_1^{2\alpha_0 +1} - \left( \frac{1}{\alpha_0} + 2\right) x_1^{\alpha_0+1} x_2^{\alpha_0}
                 + x_2^{2\alpha_0+1}}
              {(x_1^{\alpha_0} - x_2^{\alpha_0})^2}
 \end{align*}
 for \ $x_1\ne x_2$, $x_1,x_2\in I$, \ and \ $(\EE(M))(x_1,x_2) = x_1$ \ for \ $x_1=x_2\in I$.%

\noindent In case of \ $\alpha_0=0$, \ for the expectation \ $\EE(M)\in\cM_2$ \ of \ $M$, \ by \eqref{help_random_mean_calculation} and partial integration, we have
\begin{align*}
     (\EE(M))(x_1,x_2)
      = \int_0^1 x_1^\lambda x_2^{1-\lambda}  2\lambda \,\dd\lambda
       = 2x_2 \int_0^1 \lambda \left(\frac{x_1}{x_2}\right)^\lambda \,\dd\lambda
      = \frac{2(x_1\ln(x_1) - x_1 - x_1\ln(x_2) + x_2)}{(\ln(x_1) - \ln(x_2))^2}
\end{align*}
 for \ $x_1\ne x_2$, $x_1,x_2\in I\cap (0,\infty)$, \ and \ $(\EE(M))(x_1,x_2) = x_1$ \ for \ $x_1=x_2\in I\cap(0,\infty)$.
If \ $x_1=0$ \ or \ $x_2=0$, \ then \ $(\EE(M))(x_1,x_2) = 0$.%

Next, let us suppose that the distribution \ $\PP_\xi$ \ of \ $\xi$ \ takes the form
 \ $\PP_\xi = \delta_0\otimes \PP_W$, \ where \ $W$ \ is a random variable with density function
 \ $f_W(\lambda) = \frac{\ee}{\ee-1}\ee^{-\lambda} \bone_{(0,1)}(\lambda)$, \ $\lambda\in\RR$.
\ In Example \ref{Ex_2}, let us choose \ $I:=[0,b]$, \ where \ $b>0$.
\ For the expectation \ $\EE(M)\in\cM_2$ \ of \ $M$, \ by \eqref{help_random_mean_calculation}, we have
 \begin{align*}
     (\EE(M))(x_1,x_2)
      = \frac{\ee}{\ee-1} \int_0^1 x_1^\lambda x_2^{1-\lambda}  \ee^{-\lambda} \,\dd\lambda
      = \frac{\ee}{\ee-1} x_2 \int_0^1 \left(\frac{x_1}{\ee x_2}\right)^\lambda \,\dd\lambda
      = \frac{1}{\ee-1}\cdot\frac{x_1 - \ee x_2}{\ln(x_1) - \ln(\ee x_2)}
 \end{align*}
 for \ $x_1\ne \ee x_2$, \ $x_1,x_2\in I\cap (0,\infty)$, \ and \ $(\EE(M))(x_1,x_2) = \frac{1}{\ee-1}x_1$ \ for \ $x_1=\ee x_2$, $x_1,x_2\in I\cap(0,\infty)$.
If \ $x_1=0$ \ or \ $x_2=0$, \ then \ $(\EE(M))(x_1,x_2) = 0$.
\ Note that the restriction of \ $\EE(M)$ \ onto \ $(0,\infty)\times(0,\infty)$ \ can be considered as a variant of the logarithmic mean.
Namely,
 \[
 (\EE(M))(x_1,x_2)=\cL\left(\frac{x_1}{\ee-1},\frac{\ee x_2}{\ee-1}\right), \qquad  x_1,x_2\in I\cap(0,\infty) = (0,b],
 \]
 where $\cL$ denotes the logarithmic mean. However,
 note also that the mapping $I^2 \ni (x_1,x_2)\mapsto (\EE(M))(x_1,x_2)$ is a mean on its own right,
 following from Theorem \ref{Pro_random_mean_exp},  or it can be also checked directly.
Indeed, using that the function $ (0,1]\ni v\mapsto (\frac{\ee}{v}-1)/\ln(\frac{\ee}{v})$ is monotone decreasing,
 \ $ (0,1]\ni v\mapsto (v-\ee)/\ln(\frac{v}{\ee})$ \ is monotone increasing, and  that their  value at $1$ is $\ee-1$,
 we have
 \[
   1 \leq \frac{1}{\ee-1}\cdot \frac{\frac{\ee}{v}-1}{\ln\left(\frac{\ee}{v}\right)} , \qquad v\in(0,1)
    \qquad \text{and}\qquad
     \frac{1}{\ee-1} \cdot\frac{v-\ee}{\ln\left(\frac{v}{\ee}\right)}\leq 1,
     \qquad v\in(0,1),
 \]
 and hence in case of \ $0<x_1<x_2\leq b$, \ by choosing \ $v:=\frac{x_1}{x_2}$, \ we have
 \[
   x_1 \leq \frac{1}{\ee-1} \cdot \frac{x_1 - \ee x_2}{\ln(x_1) - \ln(\ee x_2)} \leq x_2,
 \]
 as desired.
Further, if \ $x_1=\ee x_2$, \ $x_1,x_2\in (0,b]$, \ then
 \[
  x_2 = \min(x_1,x_2) \leq \frac{1}{\ee-1} x_1 \leq x_1=\max(x_1,x_2),
 \]
 as desired.

Next, let us suppose that the distribution \ $\PP_\xi$ \ of \ $\xi$ \ takes the form
 \ $\PP_\xi = \delta_0\otimes \PP_X$, \ where \ $X$ \ is a random variable with density function
 \ $f_X(\lambda) = \frac{1}{1-\cos(1)}\sin(\lambda) \bone_{(0,1)}(\lambda)$, \ $\lambda\in\RR$.
\ In Example \ref{Ex_2}, let us choose \ $I:=[0,b]$, \ where \ $b>0$.
\ For the expectation \ $\EE(M)\in\cM_2$ \ of \ $M$, \ by \eqref{help_random_mean_calculation} and partial integration, we have
 \begin{align*}
     &(\EE(M))(x_1,x_2)
        = \frac{1}{1-\cos(1)} \int_0^1 x_1^\lambda x_2^{1-\lambda}  \sin(\lambda) \,\dd\lambda
        = \frac{x_2}{1-\cos(1)} \int_0^1 \left(\frac{x_1}{x_2}\right)^\lambda \sin(\lambda) \,\dd\lambda\\
     & = \frac{x_2}{1-\cos(1)} \left( 1 - \cos(1) \frac{x_1}{x_2} + \ln\left( \frac{x_1}{x_2}\right)
                                                             \int_0^1 \left(\frac{x_1}{x_2}\right)^\lambda \cos(\lambda)\,\dd \lambda  \right)\\
     & = \frac{1}{1-\cos(1)} \left( x_2 - \cos(1)x_1 + \sin(1)x_1 \ln\left( \frac{x_1}{x_2}\right)
                                     - x_2 \left(\ln\left( \frac{x_1}{x_2}\right)\right)^2  \int_0^1 \left(\frac{x_1}{x_2}\right)^\lambda \sin(\lambda)\,\dd \lambda    \right)
 \end{align*}
 for \ $x_1,x_2\in I\cap (0,\infty)$.
\ Consequently, for \ $x_1,x_2\in I\cap (0,\infty)$, \ we have
 \begin{align*}
  (\EE(M))(x_1,x_2) = \frac{1}{1 -\cos(1)} \frac{x_2 - \cos(1)x_1 + \sin(1)x_1(\ln(x_1) - \ln(x_2))}{1+ (\ln(x_1) - \ln(x_2))^2 }.
 \end{align*}
If \ $x_1=0$ \ or \ $x_2=0$, \ then \ $(\EE(M))(x_1,x_2) = 0$.%
\proofend
\end{Ex}

 Motivated by Examples \ref{Ex_2} and \ref{Ex_4}, in the next remark we initiate some possible future research directions.

\begin{Rem}\label{Rem_Future_work}
\noindent  (i) Is it possible to give a set of random variables with values in \ $\RR\times (0,1)$ \
 such that the set of expectations of the corresponding random H\"older means given in Example \ref{Ex_2} coincide with \ $\cM_p$?
If the answer is yes, then characterize such a set of  random variables.
If the answer is no, then characterize the largest subset of $\cM_p$, which can  be achieved in this way.
One can pose a similar question concerning any other random mean.

\noindent  (ii)
Moreover, given a usual (non-random) mean \ $m\in\cM_p$ \ on \ $I$, \ let us characterize (possibly under some additional assumptions)
 those \ $p\times p$ matrices \ $A$ \ with real entries
 such that the mapping \ $I^p\ni (x_1,\ldots,x_p)\mapsto m( (x_1,\ldots,x_p)A )$ \ is a (usual) mean.
 Of course, if \ $A$ \ is a $p\times p$ permutation matrix, then this property holds.
Further, in Example \ref{Ex_4}, we have showed that the mapping $(0,b]^2 \ni (x_1,x_2)\mapsto \cL((x_1,x_2)A)$ is a mean,
 where \ $b>0$, \ $\cL$ \ denotes the logarithmic mean and \ $A$ \ is the $2\times 2$ diagonal matrix with \ $(1,1)$-entry \ $\frac{1}{\ee-1}$ \ and
 \ $(2,2)$-entry \ $\frac{\ee}{\ee-1}$  (not being a permutation matrix).
\proofend
\end{Rem}

\section{Limit theorems for random means}\label{Section_Limit}

Let \ $(\Omega,\cA,\PP)$ \ be a probability space, $I$ \ be a nondegenerate, compact interval of \ $\RR$,
 \ and \ $d,p\in\NN$.
\ Let \ $(\xi_n)_{n\in\NN}$ \ be a sequence of  independent and identically distributed \ $d$-dimensional random variables,
 and for each \ $n\in\NN$, \ let \ $M_n:\Omega\to\cM_p$ \ be a random mean generated by \ $\xi_n$  (in the sense of Definition \ref{D:random_mean}).
\ For each \ $n\in\NN$, \ let \ $\overline S_n:\Omega\to\cM_p$,
 \[
   \overline S_n(\omega):= \frac{1}{n}\sum_{j=1}^n M_j(\omega),\qquad \omega\in\Omega.
 \]
Then, by Definition \ref{D:random_mean}, for each \ $n\in\NN$ \ and  \ $\omega\in\Omega$, \ we have \ $\overline S_n(\omega)\in\cM_p$, \ i.e., it is a (usual) \ $p$-variable mean,
 and for each \ $n\in\NN$, \ the mapping \ $\Omega\ni\omega\mapsto \overline S_n(\omega)$ \ is an \ $\cM_p$-valued random variable.

In what follows we are searching for sufficient conditions on the random means \ $M_n$, \ $n\in\NN$, \ under which
 \[
  \PP\Big(\big\{ \omega\in\Omega : \lim_{n\to\infty} \varrho(\overline S_n(\omega), \EE(M_1)) = 0 \big\}\Big)=1
 \]
 holds, where the  metric \ $\varrho$ \ is given in \eqref{Rem_metric}, and the law of the random variable
 \begin{align}\label{kappa}
    \Omega\ni\omega \mapsto \sqrt{n} \, \varrho(\overline S_n(\omega), \EE(M_1))=: \sqrt{n}\, \kappa_n(\omega)
 \end{align}
 converges in distribution to some normal distribution as \ $n\to\infty$.
Here note that for each \ $n\in\NN$, \ $\kappa_n$ \ is  indeed a real-valued random variable, since \ $\Omega\ni\omega\mapsto \overline S_n(\omega)$ \ is
 an \ $\cM_p$-valued random variable and the metric \ $\varrho$ \ is continuous.

\begin{Thm}\label{Pro_SLLN_CLT}
 Let \ $d,k\in\NN$, \ and \ $I$ \ be a nondegenerate, compact interval of \ $\RR$.
\ Let \ $(\xi_n)_{n\in\NN}$ \ be a sequence of independent and identically distributed \ $d$-dimensional
 discrete random variables having finite range \ $\{a_1,\ldots,a_k\}$, \ where \ $a_1,\ldots,a_k\in\RR^d$ \ are pairwise distinct.
\ Let \ $q_i:=\PP(\xi_1 = a_i)\in(0,1)$, \ $i=1,\ldots,k$.
\ Let \ $p\in \NN$ \ and for each \ $n\in\NN$ \ let \ $M_n: \Omega\to \cM_p$,
 \[
   (M_n(\omega))(x_1,\ldots,x_p):=\sum_{i=1}^k m_i(x_1,\ldots,x_p)\bone_{\{ \xi_n(\omega) = a_i\}},
    \qquad \omega\in\Omega, \;\; x_1,\ldots,x_p\in I,
 \]
 where \ $m_i\in\cM_p$, \ $i=1,\ldots,k$.
\ Then for each \ $n\in\NN$, \ $M_n$ \ is a random mean generated by \ $\xi_n$, \ and
 \begin{align}\label{help6}
  \PP\Big(\big\{ \omega\in\Omega : \lim_{n\to\infty} \kappa_n(\omega) =0 \big\}\Big)=1,
 \end{align}
 where \ $\kappa_n$ \ is given in \eqref{kappa} with \ $(\EE(M_1))(x_1,\ldots,x_p)= \sum_{i=1}^k m_i(x_1,\ldots,x_p) \PP(\xi_1=a_i)$ \
 for \ $x_1,\ldots,x_p\in I$.

Further,
 \begin{align}\label{help7}
  \sqrt{n}\,\kappa_n
  \distr \sup_{x_1,\ldots,x_p\in I}\, \big\vert \langle \cN_k(\bzero, \bQ), \bm(x_1,\ldots,x_p) \rangle \big\vert
 \end{align}
 as \ $n\to\infty$, \ where \ $\bQ := (q_{i,j})_{i,j=1}^k\in\RR^{k\times k}$ \ is the \ $k\times k$ matrix given by
 \[
   q_{i,j}:=\begin{cases}
              q_i(1-q_i) & \text{if \ $i=j$,}\\
              -q_iq_j & \text{if \ $i\ne j$,}\\
            \end{cases}
 \]
 and \ $\bm(x_1,\ldots,x_p):=( m_1(x_1,\ldots,x_p), \ldots, m_k(x_1,\ldots,x_p) )^\top\in\RR^k$.
\end{Thm}

Here \ $\bQ$ \ is nothing else but the covariance matrix of \ $(\bbone_{\{\xi_1=a_1\}},\ldots,\bbone_{\{\xi_1=a_k\}})^\top$ \ having
 multinomial distribution with parameters \ $1$ \ and \ $q_1,\ldots,q_k$.

Next, we consider a special case of  Theorem \ref{Pro_SLLN_CLT}, namely, when \ $p=2$, \ $\xi_1$ \ has a Bernoulli distribution
 and the range of \ $M_1$ \ is the set consisting of the arithmetic and geometric means in \ $[0,1]$.

\begin{Cor}\label{Cro_Bernoulli}
Let \ $I:=[0,1]$, \ and \ $(\xi_n)_{n\in\NN}$ \ be a sequence of
 independent  and identically distributed random variables such that \ $\PP(\xi_1=0)=q$ \ and \ $\PP(\xi_1=1)=1-q$, \ where \ $q\in(0,1)$,
 \  i.e., \ $\xi_1$ \ is Bernoulli distributed with parameter \ $q$.
\ For each \ $n\in\NN$, \ let \ $M_n: \Omega\to \cM_2$,
 \[
   (M_n(\omega))(x_1,x_2):=m_0(x_1,x_2)\bone_{\{ \xi_n(\omega) = 0\}} + m_1(x_1,x_2)\bone_{\{ \xi_n(\omega) = 1\}},
                           \qquad \omega\in\Omega, \;\; x_1,x_2\in I,
 \]
 where
 \[
    m_0(x_1,x_2):=\frac{x_1+x_2}{2}\qquad \text{and}\qquad m_1(x_1,x_2):=\sqrt{x_1x_2}
 \]
 for \ $x_1,x_2\in I$.
\ Then for each \ $n\in\NN$, \ $M_n$ \ is a random mean generated by \ $\xi_n$, \ and
 \begin{align}\label{help2}
  \PP\Big(\big\{ \omega\in\Omega : \lim_{n\to\infty} \kappa_n(\omega) =0 \big\}\Big)=1,
 \end{align}
 where \ $\kappa_n$ \ is given in \eqref{kappa} with \ $(\EE(M_1))(x_1,x_2)= \frac{x_1+x_2}{2} \, q + \sqrt{x_1x_2}\,(1-q)$, $x_1,x_2\in I$.
\ Further,
 \begin{align}\label{help3}
  \sqrt{n}\,\kappa_n
  \distr \sup_{x_1,x_2\in [0,1]} \left( \frac{x_1+x_2}{2} - \sqrt{x_1x_2} \right)
              \cdot \vert\cN(0,q(1-q))\vert
            = \frac{1}{2}\, \vert \cN(0,q(1-q))\vert
 \end{align}
 as \ $n\to\infty$.
\end{Cor}

Next, we establish limit theorems for randomly weighted arithmetic means, where \ $\xi_1$ \ is not necessarily discrete,
 so our next result is out of scope of  Theorem \ref{Pro_SLLN_CLT}.

\begin{Thm}\label{Pro_uniform}
Let $p\geq 2$, \ $p\in\NN$, \ $I$ \ be a nondegenerate, compact interval of \ $\RR$,
 \ and \ $(\bxi_n:=(\xi_n^{(1)},\ldots,\xi_n^{(p-1)}))_{n\in\NN}$ \ be a sequence of independent and identically distributed
 \ $\RR_+^{p-1}$-valued random variables such that \ $\PP( \xi_1^{(1)}+\cdots + \xi_1^{(p-1)} \leq 1)=1$.
\ For each \ $n\in\NN$, \ let \ $M_n:\Omega\to \cM_p$,
 \[
   (M_n(\omega))(x_1,\ldots,x_p)
    := \xi_n^{(1)}(\omega)x_1 + \cdots + \xi_n^{(p-1)}(\omega) x_{p-1} + \big(1-\xi_n^{(1)}(\omega)-\cdots - \xi_n^{(p-1)}(\omega)\big)x_p
 \]
 for $\omega\in\Omega$ \ and  $x_1,\ldots,x_p\in I$.
\ Then for each \ $n\in\NN$, \ $M_n$ \ is a random mean generated by \ $\bxi_n$, \ and
 \begin{align}\label{help12}
  \PP\Big(\big\{ \omega\in\Omega : \lim_{n\to\infty} \kappa_n(\omega) =0 \big\}\Big)=1,
 \end{align}
 where \ $\kappa_n$ \ is given in \eqref{kappa} with
 \[
   (\EE(M_1))(x_1,\ldots,x_p)
    = \sum_{i=1}^{p-1} x_i \EE(\xi_1^{(i)}) + x_p \left( 1 - \sum_{i=1}^{p-1} \EE(\xi_1^{(i)}) \right),
    \qquad x_1,\ldots,x_p\in I.
 \]
Further,
 \begin{align}\label{help13}
  \sqrt{n}\,\kappa_n
  \distr \sup_{x_1,\ldots,x_p\in I}
         \left\vert  \left\langle \cN_{p-1}\left(\bzero, \cov(\bxi_n)\right),   \begin{pmatrix}
                                                                                                       x_1-x_p \\
                                                                                                       \vdots \\
                                                                                                       x_{p-1} - x_p \\
                                                                                                     \end{pmatrix}
           \right\rangle  \right\vert
 \end{align}
 as \ $n\to\infty$, \ where \ $\cov(\bxi_n)$ \ denotes the covariance matrix of \ $\bxi_n$.
\end{Thm}

Next, we formulate a corollary of Theorem \ref{Pro_uniform} in case of \ $p=2$ \ by simplifying the limit distribution in \eqref{help13}.

\begin{Cor}\label{Cor_uniform}
Let \ $I$ \ be a nondegenerate, compact interval of \ $\RR$, \ and \ $(\xi_n)_{n\in\NN}$ \
 be a sequence of independent and identically distributed random variables such that \ $\PP(\xi_1\in [0,1])=1$.
\ For each \ $n\in\NN$, \ let \ $M_n:\Omega\to \cM_2$, \ $(M_n(\omega))(x_1,x_2):=\xi_n(\omega)x_1 + (1-\xi_n(\omega))x_2$,
 \ $\omega\in\Omega$, \ $x_1,x_2\in I$.
\ Then for each \ $n\in\NN$, \ $M_n$ \ is a random mean generated by \ $\xi_n$, \ and
 \begin{align}\label{help4}
  \PP\Big(\big\{ \omega\in\Omega : \lim_{n\to\infty} \kappa_n(\omega) =0 \big\}\Big)=1,
 \end{align}
 where \ $\kappa_n$ \ is given in \eqref{kappa} with \ $(\EE(M_1))(x_1,x_2) = x_1\EE(\xi_1) + x_2(1-\EE(\xi_1))$, \ $x_1,x_2\in I$.
\ Further,
 \begin{align}\label{help5}
  \sqrt{n}\,\kappa_n
  \distr \sup_{x_1,x_2\in I} \vert x_1-x_2\vert\cdot
         \left\vert \cN\left(0,\DD^2(\xi_1)\right)\right\vert
 \end{align}
 as \ $n\to\infty$.
\end{Cor}

Finally, we provide limit theorems for randomly weighted power means (which can be also called random H\"older means, see Example \ref{Ex_2}).
We point out to the facts that in this case instead of the arithmetic mean of the given random means we consider their geometric mean, and so
 the limit theorems have somewhat different forms compared to the previous ones.

\begin{Thm}\label{Pro_weighted_power_mean}
Let \ $I$ \ be a nondegenerate, compact interval  of \ $(0,\infty)$, \ and
 \ $(\xi_n)_{n\in\NN}$ \ be a sequence of independent and identically distributed random variables such that \ $\PP(\xi_1\in[0,1])=1$.
\ For each \ $n\in\NN$, \ let \ $M_n:\Omega\to \cM_2$, \ $(M_n(\omega))(x_1,x_2):=x_1^{\xi_n(\omega)}x_2^{1-\xi_n(\omega)}$, \ $\omega\in\Omega$, \ $x_1,x_2\in I$.
\ Then for each \ $n\in\NN$, \ $M_n$ \ is a random mean generated by \ $\xi_n$, \ and
 \begin{align}\label{help10}
  \PP\left(\left\{ \omega\in\Omega : \lim_{n\to\infty} \sup_{x_1,x_2\in I}
                \frac{\left( \prod_{j=1}^n (M_j(\omega))(x_1,x_2) \right)^{\frac{1}{n}}}
                      {x_1^{\EE(\xi_1)} x_2^{1-\EE(\xi_1)}}
                  = 1 \right\}\right)=1.
 \end{align}
Further,
 \begin{align}\label{help11}
  \left( \sup_{x_1,x_2\in I}
         \frac{\left( \prod_{j=1}^n (M_j(\cdot))(x_1,x_2)\right)^{\frac{1}{n}}}
                      {x_1^{\EE(\xi_1)} x_2^{1-\EE(\xi_1)}}
  \right)^{\sqrt{n}}
  \distr \left(\frac{ \max(I)}{ \min(I)}\right)^{\left\vert \cN\left(0,\DD^2(\xi_1)\right)\right\vert}
 \end{align}
 as \ $n\to\infty$, \ where for any \ $j\in\NN$ \ and \ $x_1,x_2\in I$, \ $(M_j(\cdot))(x_1,x_2)$ \ denotes the random variable
 \ $\Omega\ni\omega\mapsto (M_j(\omega))(x_1,x_2)$, \ and \ $\max(I):=\max\{x : x\in I\}$ \ and \ $\min(I):=\min\{x : x\in I\}$.
\end{Thm}

\begin{Rem}\label{Rem_Frechet}
{Using the notations of the second part of Remark \ref{Rem_random_mean_comparison}, we} note that
 a Wasserstein law of large numbers holds for a sequence \ $(\Lambda_n)_{n\in\NN}$ \ of independent and identically
 distributed random means with values in \ $\cW_2(\RR^d)$ \ having unique Fr\'echet means, namely, the so-called empirical Fr\'echet mean
 of \ $(\Lambda_1,\ldots,\Lambda_n)$ \ (see Panaretos and Zemel \cite[Definition 3.1.1]{PanZem}) converges almost surely to the Fr\'echet mean
 of \ $\Lambda_1$ \ as \ $n\to\infty$, \ see Corollary 3.2.10 in Panaretos and Zemel \cite{PanZem}.
Note that in present section, we derived different kinds of limit theorems for random means generated by a sequence of independent and identically
 distributed random variables, since our limit theorems are about the random means itself and not about their expectations.
\proofend
\end{Rem}

\section{Proofs for Section \ref{Section_Limit}}\label{Section_proofs}

\noindent{\bf Proof of Theorem \ref{Pro_SLLN_CLT}.}
By Example \ref{Ex_1}, for each \ $n\in\NN$, \ $M_n$ \ is a random mean generated by \ $\xi_n$.
\ For all \ $x_1,\ldots,x_p\in I$, \  we have
 \begin{align*}
   (\overline S_n(\omega))(x_1,\ldots,x_p)
    &=\frac{1}{n}\sum_{j=1}^n \bbone_{\{\xi_j(\omega) = a_1\}}\, m_1(x_1,\ldots,x_p)
      + \cdots + \frac{1}{n}\sum_{j=1}^n \bbone_{\{\xi_j(\omega) = a_k\}} \, m_k(x_1,\ldots,x_p),
 \end{align*}
 and, by Example \ref{Ex_3},
 \[
   (\EE(M_1))(x_1,\ldots,x_p) = \sum_{i=1}^k m_i(x_1,\ldots,x_p) \PP(\xi_1=a_i),
   \qquad x_1,\ldots,x_p\in I.
 \]
Hence
 \begin{align*}
  &\kappa_n(\omega) = \varrho(\overline S_n(\omega), \EE(M_1))
     = \sup_{(x_1,\ldots,x_p)\in I^p} \left \vert  (\overline S_n(\omega))(x_1,\ldots,x_p) -  (\EE(M_1))(x_1,\ldots,x_p) \right\vert =
  \end{align*}
 \begin{align*}
  &= \sup_{(x_1,\ldots,x_p)\in I^p} \Bigg \vert  \left( \frac{1}{n}\sum_{j=1}^n \bbone_{\{\xi_j(\omega) = a_1\}} - \PP(\xi_1=a_1) \right) m_1(x_1,\ldots,x_p) \\
  &\phantom{= \sup_{(x_1,\ldots,x_p)\in I^p} \quad} + \cdots + \left( \frac{1}{n}\sum_{j=1}^n \bbone_{\{\xi_j(\omega) = a_k\}} - \PP(\xi_1=a_k) \right) m_k(x_1,\ldots,x_p)
             \Bigg\vert \\
  &\leq \left\vert \frac{1}{n}\sum_{j=1}^n \bbone_{\{\xi_j(\omega) = a_1\}} - \PP(\xi_1=a_1) \right\vert
       \sup_{(x_1,\ldots,x_p)\in I^p} \vert m_1(x_1,\ldots,x_p) \vert \\
  & \phantom{\quad} + \cdots +
    \left\vert \frac{1}{n}\sum_{j=1}^n \bbone_{\{\xi_j(\omega) = a_k\}} - \PP(\xi_1=a_k) \right\vert
       \sup_{(x_1,\ldots,x_p)\in I^p} \vert m_k(x_1,\ldots,x_p) \vert \\
  &\to 0 \qquad \text{as \ $n\to\infty$}
 \end{align*}
 for \ $\PP$-a.e. \ $\omega\in\Omega$, \ yielding \eqref{help6}, where the last step follows by the strong law of large numbers and
 \[
   \sup_{(x_1,\ldots,x_p)\in I^p} \vert m_i(x_1,\ldots,x_p) \vert = \max_{t\in I}\vert t\vert<\infty,\qquad i=1,\ldots,k,
 \]
 where we used that \ $m_i(x_1,\ldots,x_p) = x_1$ \ whenever \ $x_1=\ldots=x_p\in I$ \ and that \ $I$ \ is compact.

Now we turn to prove \eqref{help7}.
For all \ $\omega\in\Omega$ \ and \ $x_1,\ldots,x_p\in I$, \ we have
 \begin{align*}
   &\sqrt{n} \left \vert  (\overline S_n(\omega))(x_1,\ldots,x_p) -  (\EE(M_1))(x_1,\ldots,x_p) \right\vert \\
   &\qquad  = \left\vert \left\langle  \sqrt{n} \begin{pmatrix}
                                            \frac{1}{n}\sum_{j=1}^n \bbone_{\{\xi_j(\omega) = a_1\}} - \PP(\xi_1=a_1)  \\
                                           \vdots \\
                                           \frac{1}{n}\sum_{j=1}^n \bbone_{\{\xi_j(\omega) = a_k\}} - \PP(\xi_1=a_k)  \\
                                         \end{pmatrix},
                                         \bm(x_1,\ldots,x_p)
       \right\rangle \right\vert,
 \end{align*}
 and, by the multidimensional central limit theorem,
 \begin{align*}
  \sqrt{n} \begin{pmatrix}
    \frac{1}{n}\sum_{j=1}^n \bbone_{\{\xi_j = a_1\}} - \PP(\xi_1=a_1)  \\
    \vdots \\
    \frac{1}{n}\sum_{j=1}^n \bbone_{\{\xi_j = a_k\}} - \PP(\xi_1=a_k)  \\
  \end{pmatrix}
  \distr \cN_k\Big(\bzero, (\cov(\bbone_{\{\xi_1=a_i\}}, \bbone_{\{\xi_1=a_j\}} ))_{i,j=1}^k\Big)
 \end{align*}
 as \ $n\to\infty$, \ where
 \begin{align*}
  \cov(\bbone_{\{\xi_1=a_i\}}, \bbone_{\{\xi_1=a_j\}} )
   &= \begin{cases}
        \PP(\xi_1=a_i) - \PP(\xi_1=a_i)^2 & \text{if \ $i=j$,}\\
         - \PP(\xi_1=a_i)\PP(\xi_1=a_j) & \text{if \ $i\ne j$,}\\
     \end{cases}\\
   &= q_{i,j},
   \qquad i,j=1,\ldots,k.
 \end{align*}
Further, since \ $I$ \ is compact and \ $m_i\in\cM_p$, \ $i=1,\ldots,k$,
 \ we have the set \ $\bm(I\times \cdots\times I) = I^k$ \ is a compact subset of \ $\RR^k$, \ so, by Theorem \ref{Thm_continuity},
 the mapping
 \[
 \RR^k\ni\bq \mapsto
      \sup_{\by \in \bm(I,\ldots,I) }  \vert  \langle \bq, \by \rangle  \vert
       = \sup_{(x_1,\ldots,x_p)\in I^p}  \vert  \langle \bq, \bm(x_1,\ldots,x_p)\rangle  \vert
 \]
 is well-defined and continuous.
Consequently, the continuous mapping theorem yields \eqref{help7}.
\proofend

\noindent{\bf  First proof of Corollary \ref{Cro_Bernoulli}.}
We can apply  Theorem \ref{Pro_SLLN_CLT} with $p=2$.
Namely, using the notations of  Theorem \ref{Pro_SLLN_CLT}, we have
 \[
  \cN_2(\bzero, \bQ) =  \cN_2\left( \begin{pmatrix}
                                      0 \\
                                      0 \\
                                    \end{pmatrix}
                                , \begin{pmatrix}
                                    q(1-q) & -q(1-q) \\
                                    -q(1-q) & q(1-q) \\
                                  \end{pmatrix}
                                 \right)
                       \distre     \begin{pmatrix}
                                      -\eta \\
                                      \eta \\
                                  \end{pmatrix},
 \]
 where \ $\eta$ \ is a 1-dimensional random variable having distribution \ $\cN(0,q(1-q))$, \ and
 \[
   \bm(x_1,x_2) = \begin{pmatrix}
                    \frac{x_1+x_2}{2} \\
                    \sqrt{x_1x_2} \\
                  \end{pmatrix},
                  \qquad x_1,x_2\in [0,1].
 \]
Hence
 \begin{align*}
  &\sup_{x_1,x_2\in [0,1]}\, \big\vert \langle \cN_2(\bzero, \bQ), \bm(x_1,x_2) \rangle \big\vert
   \distre \sup_{x_1,x_2\in [0,1]} \left\vert -\eta \frac{x_1+x_2}{2} + \eta \sqrt{x_1x_2} \right\vert \\
  &\qquad  = \vert \eta\vert \sup_{x_1,x_2\in [0,1]} \left \vert \frac{x_1+x_2}{2} - \sqrt{x_1x_2}\right \vert
    \distre  \frac{1}{2} \vert \cN(0,q(1-q))\vert,
 \end{align*}
 as desired, since
 \begin{align}\label{help9}
  \begin{split}
    &\sup_{x_1,x_2\in [0,1]}   \left\vert \frac{x_1+x_2}{2} - \sqrt{x_1x_2} \right\vert
          = \sup_{x_1,x_2\in [0,1]}   \left(\frac{x_1+x_2}{2} - \sqrt{x_1x_2} \right)  \\
    & \qquad = \frac{1}{2} \sup_{x_1,x_2\in [0,1]} (\sqrt{x_1} - \sqrt{x_2})^2
             = \frac{1}{2} \left(  \sup_{x_1,x_2\in [0,1]} \vert \sqrt{x_1} - \sqrt{x_2} \vert \right)^2
             = \frac{1}{2}.
  \end{split}
 \end{align}
\proofend

\noindent{\bf Second proof of Corollary \ref{Cro_Bernoulli}.}
We give a direct proof as well, not refereeing to  Theorem \ref{Pro_SLLN_CLT}.
By Example \ref{Ex_1}, for each \ $n\in\NN$, \ $M_n$ \ is a random mean generated by \ $\xi_n$.
\  Using that \ $\sum_{j=1}^n \bbone_{\{\xi_j = 1\}}  = n - \sum_{j=1}^n \bbone_{\{\xi_j = 0\}}$,
 \ $n\in\NN$, \ for all \ $x_1,x_2\in I$, \ we have
 \begin{align*}
   (\overline S_n(\omega))(x_1,x_2)
    =\frac{1}{n}\sum_{j=1}^n \bbone_{\{\xi_j(\omega) = 0\}} \frac{x_1+x_2}{2}
      + \frac{1}{n}\sum_{j=1}^n \bbone_{\{\xi_j(\omega) = 1\}} \sqrt{x_1x_2}\\
    = \frac{1}{n} \sum_{j=1}^n \bbone_{\{\xi_j(\omega) = 0\}} \left( \frac{x_1+x_2}{2} - \sqrt{x_1x_2} \right)
      + \sqrt{x_1x_2},
 \end{align*}
 and, by Example \ref{Ex_3},
 \begin{align*}
   (\EE(M_1))(x_1,x_2) & = \frac{x_1+x_2}{2} \PP(\xi_1=0)  + \sqrt{x_1x_2} \PP(\xi_1=1) \\
                       & = \PP(\xi_1 = 0) \left( \frac{x_1+x_2}{2} - \sqrt{x_1x_2} \right) + \sqrt{x_1x_2},
                           \qquad x_1,x_2\in I.
 \end{align*}
Hence
 \begin{align*}
   \kappa_n(\omega)
     & = \varrho(\overline S_n(\omega), \EE(M_1) )
        = \sup_{(x_1,x_2)\in I^2} \left \vert  (\overline S_n(\omega))(x_1,x_2) -  (\EE(M_1))(x_1,x_2) \right\vert \\
     & = \sup_{(x_1,x_2)\in I^2}  \left\vert \left(  \frac{1}{n}\sum_{j=1}^n \bbone_{\{\xi_j(\omega) = 0\}}  - \PP(\xi_1=0) \right)
                                 \left( \frac{x_1+x_2}{2} - \sqrt{x_1x_2} \right)  \right\vert \\
     & = \left(\sup_{(x_1,x_2)\in I^2}   \left( \frac{x_1+x_2}{2} - \sqrt{x_1x_2} \right) \right)
              \left\vert \frac{1}{n}\sum_{j=1}^n \bbone_{\{\xi_j(\omega) = 0\}}  - \PP(\xi_1=0) \right\vert
 \end{align*}
 for all  $n\in\NN$ \ and\ $\omega\in\Omega$.
\ By the strong law of large numbers, we have \eqref{help2}.
The central limit theorem together with \eqref{help9} and the continuous mapping theorem applied to the function \ $\RR\ni x\mapsto \vert x\vert$ \
 yield \eqref{help3}.
\proofend

\noindent{\bf Proof of Theorem \ref{Pro_uniform}.}
First, we check that for each \ $n\in\NN$, \ $M_n$ \ is a random mean generated by \ $\bxi_n$.
\ For each \ $n\in\NN$ \ and \ $\omega\in\Omega$, \ $M_n(\omega)$ \ can be written in the form
 \[
   (M_n(\omega))(x_1,\ldots,x_p) = f(x_1,\ldots,x_p,\bxi_n(\omega)), \qquad  x_1,\ldots,x_p\in I,
 \]
 where \ $f:I^p\times \RR^{p-1}\to I$ \ is a
 \ $(\cB(I^p)\times \cB(\RR^{p-1}),\cB(I))$-measurable function satisfying
 \[
    f(x_1,\ldots,x_p,y_1,\ldots,y_{p-1}) = y_1x_1 + \cdots + y_{p-1} x_{p-1} + (1-y_1-\cdots - y_{p-1})x_p
 \]
 for \ $x_1,\ldots,x_p\in I$, \ $y_1,\ldots, y_{p-1}\in\RR_+$ \ with \ $y_1+\cdots+y_{p-1}\leq 1$,
 \ and \ $f(\cdot,\ldots,\cdot,y_1,\ldots,y_{p-1})$ \ is a fixed (arbitrary) element of \ $\cM_p$ \
 for any \ $(y_1,\ldots,y_{p-1}) \in \RR^{p-1}\setminus\{ (y_1,\ldots,y_{p-1})\in\RR_+^p : y_1+\cdots+y_{p-1}\leq 1\}$.
\ Hence, by Theorem \ref{Pro_3}, \ $M_n$ \ is a random mean generated by \ $\bxi_n$ \ for each \ $n\in\NN$.
\ Further, for the expectation \ $\EE(M_1)\in\cM_p$ \ of \ $M_1$ \ we have
 \begin{align*}
   & (\EE(M_1))(x_1,\ldots,x_p)  = \int_{\RR^{p-1}} f(x_1,\ldots,x_p,y_1,\ldots,y_{p-1})\,\PP_{\bxi_1}(\dd y_1,\ldots,\dd y_{p-1}) \\
   & = \int_{\{ y_1,\ldots,y_{p-1}\in\RR_+ : y_1+\cdots+y_{p-1}\leq 1 \}  }
           \Big(y_1x_1+\cdots + y_{p-1}x_{p-1} + (1-y_1-\cdots -y_{p-1})x_p \Big)\,\PP_{\bxi_1}(\dd y_1,\ldots,\dd y_{p-1})\\
   & = \sum_{i=1}^{p-1} x_i \EE(\xi_1^{(i)}) + x_p \left( 1 - \sum_{i=1}^{p-1} \EE(\xi_1^{(i)}) \right),
                               \qquad x_1,\ldots,x_p\in I,
 \end{align*}
 where \ $ \EE(\xi_1^{(i)})\in[0,1]$, \ $i=1,\ldots,p-1$, \  $1 - \sum_{i=1}^{p-1} \EE(\xi_1^{(i)}) \in [0,1]$, \
 and the first equality follows by \eqref{help_random_mean_calculation}.

For all  \ $n\in\NN$, \ $\omega\in\Omega$ \ and \ $x_1,\ldots,x_p\in I$, \ we have
 \begin{align*}
   (\overline S_n(\omega))(x_1,\ldots,x_p)
    & = \frac{1}{n} \sum_{j=1}^n \Big( \xi^{(1)}_j(\omega) x_1 + \cdots + \xi^{(p-1)}_j(\omega) x_{p-1} + (1 - \xi^{(1)}_j(\omega) - \cdots - \xi^{(p-1)}_j(\omega)) x_p \Big) =
 \end{align*}
 \begin{align*}
    & = \left(\frac{1}{n} \sum_{j=1}^n \xi^{(1)}_j(\omega) \right) x_1 + \cdots + \left(\frac{1}{n} \sum_{j=1}^n \xi^{(p-1)}_j(\omega)\right) x_{p-1}\\
    &\phantom{=\;}   + \left(\! 1 - \frac{1}{n} \sum_{j=1}^n \xi^{(1)}_j(\omega) - \cdots - \frac{1}{n} \sum_{j=1}^n \xi^{(p-1)}_j(\omega) \right) x_p.
 \end{align*}
Hence
 \begin{align*}
   &\kappa_n(\omega)=\varrho(\overline S_n(\omega), \EE(M_1) )\\
   &\qquad = \sup_{x_1,\ldots,x_p\in I} \left\vert  (\overline S_n(\omega))(x_1,\ldots,x_p)
                                             -  (\EE(M_1))(x_1,\ldots,x_p)   \right\vert \\
   &\qquad = \sup_{x_1,\ldots,x_p\in I}
       \Bigg\vert
         \left(  \frac{1}{n} \sum_{j=1}^n \xi^{(1)}_j(\omega) - \EE(\xi^{(1)}_1) \right) x_1
                 + \cdots + \left(  \frac{1}{n} \sum_{j=1}^n \xi^{(p-1)}_j(\omega) - \EE(\xi^{(p-1)}_1) \right) x_{p-1}\\
   &\phantom{\qquad = \sup_{x_1,\ldots,x_p\in I} \Bigg\vert \Bigg(}
                 - \left[ \frac{1}{n} \sum_{j=1}^n \xi^{(1)}_j(\omega) - \EE(\xi^{(1)}_1) + \cdots + \frac{1}{n} \sum_{j=1}^n \xi^{(p-1)}_j(\omega) - \EE(\xi^{(p-1)}_1) \right]x_p
       \,\Bigg\vert \\
   &\qquad \leq 2 \left( \max_{t\in I} \vert t\vert\right) \sum_{k=1}^{p-1} \left\vert \frac{1}{n} \sum_{j=1}^n \xi^{(k)}_j(\omega)  - \EE(\xi^{(k)}_1)  \right\vert
     \to 0 \qquad \text{as \ $n\to\infty$}
 \end{align*}
 for \ $\PP$-a.e. \ $\omega\in\Omega$, \ yielding \eqref{help12}, where we used the strong law of large numbers and that
 \ $\max_{t\in I} \vert t\vert <\infty$, \ since \ $I$ \ is compact.

Now we turn to prove \eqref{help13}.
For all \ $\omega\in\Omega$ \ and \ $x_1,\ldots,x_p\in I$, \ we have
 \begin{align*}
  &\sqrt{n} \Big\vert (\overline S_n(\omega))(x_1,\ldots,x_p)
                 -  (\EE(M_1))(x_1,\ldots,x_p) \Big \vert \\
  &\qquad
    = \left\vert
    \left\langle
     \sqrt{n}
     \begin{pmatrix}
        \frac{1}{n} \sum_{j=1}^n \xi^{(1)}_j(\omega) - \EE(\xi^{(1)}_1) \\
        \vdots \\
       \frac{1}{n} \sum_{j=1}^n \xi^{(p-1)}_j(\omega) - \EE(\xi^{(p-1)}_1) \\
    \end{pmatrix}  ,
    \begin{pmatrix}
       x_1-x_p \\
        \vdots \\
       x_{p-1} - x_p \\
    \end{pmatrix}
  \right\rangle
  \right\vert,
 \end{align*}
 and, by the multidimensional central limit theorem,
 \[
   \sqrt{n}
     \begin{pmatrix}
        \frac{1}{n} \sum_{j=1}^n \xi^{(1)}_j - \EE(\xi^{(1)}_1) \\
        \vdots \\
       \frac{1}{n} \sum_{j=1}^n \xi^{(p-1)}_j - \EE(\xi^{(p-1)}_1) \\
    \end{pmatrix}
    \distr
    \cN_{p-1} (\bzero, \cov(\xi^{(1)}_1,\ldots, \xi^{(p-1)}_1))
 \]
 as \ $n\to\infty$.
Since \ $I$ \ is compact, the set
 \[
   \left\{   \begin{pmatrix}
       x_1-x_p \\
        \vdots \\
       x_{p-1} - x_p \\
    \end{pmatrix}
      : x_1,\ldots,x_p\in I  \right\}
 \]
 is compact as well.
Indeed,
 \[
   \begin{pmatrix}
       x_1-x_p \\
        \vdots \\
       x_{p-1} - x_p \\
    \end{pmatrix}
     = \begin{pmatrix}
         1 & 0 & \cdots & 0 & -1 \\
         0 & 1 & \cdots & 0 & -1 \\
         \vdots & \vdots & \ddots & \vdots & \vdots \\
         0 & 0 & \cdots & 1 & -1 \\
       \end{pmatrix}
       \begin{pmatrix}
       x_1 \\
       x_2 \\
        \vdots \\
       x_p \\
    \end{pmatrix},
       \qquad (x_1,\ldots,x_p)\in I^p,
 \]
 the set \ $I^p$ \ is compact, and it is known that the image of a compact set of \ $\RR^p$ \ by a continuous map is compact.
Hence, by Theorem \ref{Thm_continuity}, the mapping
 \[
 \RR^{p-1}\ni \bq \mapsto \sup_{x_1,\ldots,x_p\in I}
  \left\vert \left\langle
      \bq,
     \begin{pmatrix}
       x_1-x_p \\
        \vdots \\
       x_{p-1} - x_p \\
    \end{pmatrix}
    \right\rangle
    \right\vert
 \]
 is well-defined and continuous, so the continuous mapping theorem yields \eqref{help13}.
\proofend

\noindent{\bf Proof of Corollary \ref{Cor_uniform}.}
Theorem \ref{Pro_uniform} yields Corollary \ref{Cor_uniform}, since
 \[
    \sup_{x_1,x_2\in I} \Big\vert \langle \cN(0,\DD^2(\xi_1)), x_1 - x_2  \rangle \Big\vert
    =  \sup_{x_1,x_2\in I} \vert x_1 -x_2 \vert \cdot \vert \cN(0,\DD^2(\xi_1)) \vert.
 \]
\proofend

\noindent{\bf Proof of Theorem \ref{Pro_weighted_power_mean}.}
First, we check that for each \ $n\in\NN$, \ $M_n$ \ is a random mean generated by \ $\xi_n$.
\ For all \ $n\in\NN$ \ and \ $\omega\in\Omega$, \ $M_n(\omega)$ \ can be written in the form
 \[
   (M_n(\omega))(x_1,x_2) = f(x_1,x_2,\xi_n(\omega)), \qquad  x_1,x_2\in I,
 \]
 where \ $f:I^2\times \RR\to I$ \ is a \ $(\cB(I^2)\times\cB(\RR),\cB(I))$-measurable function satisfying
 \[
    f(x_1,x_2,y)= x_1^yx_2^{1-y},\qquad x_1,x_2\in I, \;\; y\in[0,1],
 \]
 and \ $f(\cdot,\cdot,y)$ \ is a fixed (arbitrary) element of \ $\cM_2$ \
 for any \ $y\in \RR\setminus [0,1]$.
\ Hence, by Theorem \ref{Pro_3}, \ $M_n$ \ is a random mean generated by \ $\xi_n$ \ for each \ $n\in\NN$.
\ We also have that for each \ $n\in\NN$ \ and \ $x_1,x_2\in I$, \ the mapping \ $\Omega\ni\omega\mapsto (M_n(\omega))(x_1,x_2)$ \
 is \ $(\cA,\cB(\RR))$-measurable, i.e., it is a random variable, since
 \ $(M_n(\omega))(x_1,x_2) = \varphi_{x_1,x_2}(M_n(\omega))$, $\omega\in\Omega$,
 \ $M_n$ \ is \ $(\cA,\cB(\cC(I^2)))$-measurable and
 \ $\varphi_{x_1,x_2}(h):=h(x_1,x_2)$, \ $h\in\cC(I^2)$, \ is a linear functional.

By the assumptions, there exists \ $c,\widetilde c\in(0,\infty)$ \ such that \ $c< x < \widetilde c$ \ for all \ $x\in I$, \
 so \ $0< c \leq  \min(I) \leq \max(I) \leq \widetilde c <\infty$.
\ Further, note that \ $\EE(\xi_1)$ \ and \ $\DD^2(\xi_1)$ \ exist, and \ $\EE(\xi_1)\in[c,\widetilde c]$.
\ For all \ $n\in\NN$ \ and \ $x_1,x_2\in I$, \ we have
 \begin{align*}
  \frac{\left( \prod_{j=1}^n  (M_j(\cdot))(x_1,x_2) \right)^{\frac{1}{n}}}
       {x_1^{\EE(\xi_1)} x_2^{1-\EE(\xi_1)}}
  & = \frac{\left( \prod_{j=1}^n x_1^{\xi_j} x_2^{1-\xi_j} \right)^{\frac{1}{n}}}
       {x_1^{\EE(\xi_1)} x_2^{1-\EE(\xi_1)}}
     = x_1^{\frac{1}{n} \sum_{j=1}^n \xi_j - \EE(\xi_1)}
      x_2^{\frac{1}{n} \sum_{j=1}^n (1-\xi_j) - (1-\EE(\xi_1))}.
 \end{align*}
Hence, using that the functions \ $\exp$ \ and \ $\ln$ \ are strictly increasing, we have for each \ $n\in\NN$,
 \begin{align*}
  & \sup_{x_1,x_2\in I}
       \frac{\left( \prod_{j=1}^n  (M_j(\cdot))(x_1,x_2) \right)^{\frac{1}{n}}}
            {x_1^{\EE(\xi_1)} x_2^{1-\EE(\xi_1)}}
   = \sup_{x_1,x_2\in I}
       x_1^{\frac{1}{n} \sum_{j=1}^n \xi_j - \EE(\xi_1)}
       x_2^{\frac{1}{n} \sum_{j=1}^n (1-\xi_j) - (1-\EE(\xi_1))}\\
  &\qquad  = \sup_{x_1,x_2\in I}
      \exp\left\{  \left( \frac{1}{n} \sum_{j=1}^n \xi_j - \EE(\xi_1) \right)\ln(x_1)
                  + \left( \frac{1}{n} \sum_{j=1}^n (1-\xi_j) - (1-\EE(\xi_1)) \right) \ln(x_2)
          \right\}\\
  &\qquad  =  \exp\left\{  \sup_{x_1,x_2\in I} \left(  \frac{1}{n} \sum_{j=1}^n \xi_j - \EE(\xi_1) \right)   \ln\left( \frac{x_1}{x_2}\right) \right\}\\
  & \qquad =  \exp\Bigg\{
                    \left(  \frac{1}{n} \sum_{j=1}^n \xi_j - \EE(\xi_1) \right) \bbone_{\big\{  \frac{1}{n} \sum_{j=1}^n \xi_j - \EE(\xi_1) \geq 0 \big\}}
                     \sup_{x_1,x_2\in I}  \ln\left( \frac{x_1}{x_2}\right) \\
  &\phantom{\qquad =  \exp\Bigg\{ }
             + \left(  \frac{1}{n} \sum_{j=1}^n \xi_j - \EE(\xi_1) \right) \bbone_{\big\{  \frac{1}{n} \sum_{j=1}^n \xi_j - \EE(\xi_1) < 0 \big\}}
               \inf_{x_1,x_2\in I}  \ln\left( \frac{x_1}{x_2}\right)
                \Bigg\} \\
 &\qquad = \exp\left\{  \left\vert \frac{1}{n} \sum_{j=1}^n \xi_j - \EE(\xi_1) \right\vert \ln\left(  \frac{\max(I)}{\min(I)} \right) \right\}\\
 &\qquad = \left(  \frac{\max(I)}{\min(I)}\right)^{\left\vert \frac{1}{n} \sum_{j=1}^n \xi_j - \EE(\xi_1) \right\vert}.
 \end{align*}
Using the strong law of large numbers and that \ $\frac{\max(I)}{\min(I)}\in(0,\infty)$, \ we have \eqref{help10}.
The central limit theorem together with the continuous mapping theorem  applied to the function
 \ $\RR\ni x \mapsto \left(  \frac{\max(I)}{\min(I)} \right)^{\vert x\vert}$ \ yield \eqref{help11}.
\proofend

\section{Declarations}

\noindent {\bf Funding.}
M\'aty\'as Barczy is supported by grant NKFIH-1279-2/2020 of the Ministry for Innovation and Technology, Hungary.

\noindent  {\bf Conflicts of interest/Competing interests.} We do not have any conflicts of interest/competing interests.

\noindent  {\bf Availability of data and material.} Not applicable.

\noindent  {\bf Code availability.} Not applicable.

\appendix

\vspace*{5mm}

\noindent{\bf\Large Appendix}

\section{Continuity of the supremum}

The following result is known, however, we could not address any book or article containing it,
 only an internet blog due to Wong \cite{Won}.
Because it is used in the  verifications of  Theorems \ref{Pro_SLLN_CLT} and \ref{Pro_uniform} we present a proof of it.

\begin{Thm}\label{Thm_continuity}
Let \ $X$ \ and \ $Y$ \ be topological spaces such that \ $Y$ \ is compact, and let \ $f:X\times Y\to\RR$ \ be a continuous function.
Then the function \ $g:X\to\RR$, \ $g(x):=\sup\limits_{y\in Y} f(x,y)$, \ $x\in X$, \ is well-defined and continuous.
\end{Thm}

\noindent{\bf Proof.}
Because of the continuity of $f$,  the function $Y\ni y\mapsto f(x,y)$ is continuous for every fixed $x\in X$,
  and the compactness of $Y$ implies  that its supremum is finite and it is  attained.
So, $g(x)$, $ x\in X$, is well-defined.

Let $r\in\RR$ be arbitrarily fixed.
We will prove that the inverse images $g^{-1}((-\infty,r))$ and $g^{-1}((r,\infty))$ are open in $X$.

 First, we prove that $g^{-1}((r,\infty))$ is open in $X$.
\ Let us denote by $\pi_X\colon X\times Y\to X,\ \pi_X(x,y) :=x$, $ (x,y)\in X\times Y$, the canonical projection onto $X$,
 which is known to be continuous and open  (i.e., maps open sets to open sets).
Moreover, for every $x_0\in X$ there is at least one $y_0\in Y$ such that
\[
g(x_0)=\sup\limits_{y\in Y} f(x_0,y)=f(x_0,y_0).
\]
So, for every $r\in\RR$ we can write
\begin{align*}
 g^{-1}((r,\infty))& =\{x\in X\ |\ g(x)>r \}=\Big\{ x\in X\ |\ \sup\limits_{y\in Y}f(x,y)>r \Big\}\\
                   & = \big\{x\in X\ |\ f(x,y)>r  \mbox{ for some }y\in Y\big\}=  \pi_X (f^{-1}((r,\infty))).
\end{align*}
Because $f$ is continuous, $f^{-1}((r,\infty))$ is open in \ $X\times Y$.
\ The canonical projection \ $ \pi_X$ \ is an open map, which entails that $ \pi_X (f^{-1}((r,\infty)))$ is an open subset of $X$,
 so is $g^{-1}((r,\infty))$.

Next, we prove that $g^{-1}((-\infty,r))$ is open  in $X$ for every $r\in\RR$.
If $g(x)<r$  for some $x\in X$, then, by the definition of $g$, we have $f(x,y) <\tilde r<r$ for every $y\in Y$, \  where \ $\tilde r$ \ satisfies \ $g(x)<\tilde r<r$.
\ In other words, if $x\in g^{-1}((-\infty,r))$, then $\{x\}\times Y\subset f^{-1}((-\infty,\tilde r))$.
\ Because of the continuity of $f$, the set $f^{-1}((-\infty, \tilde r))$ is open in $X\times Y$.
So,  if $x\in g^{-1}((-\infty,r))$ and $y\in Y$, then there are open sets $U_{x,y}\subset X$ and $V_{x,y}\subset Y$
 such that $U_{x,y}\times V_{x,y}$ is an open neighbourhood of $(x,y)\in X\times Y$  and it is contained in $f^{-1}((-\infty, \tilde r))$.
\ For a fixed $x\in g^{-1}((-\infty,r))$, the sets $V_{x,y}$, $ y\in Y$, give an open cover of $Y$,
  and, because of the compactness of \ $Y$, there  exist $k(x)\in\NN$  and $y_1,\ldots,y_{k(x)}\in Y$
 such that $Y =\bigcup_{i=1}^{k(x)}V_{x,y_i}$.
 Using that \ $A\times (B\cup C) = (A\times B)\cup (A\times C)$ \ for any sets \ $A,B,C$, \ this entails that
 \begin{align*}
  \{x\}\times Y\subset\left(\bigcap_{i=1}^{k(x)}U_{x,y_i}\right)\times Y
                = \bigcup_{j=1}^{k(x)}  \left(\left(\bigcap_{i=1}^{k(x)}U_{x,y_i} \right) \times V_{x,y_j} \right)
                     \subset \bigcup_{j=1}^{k(x)} \left(U_{x,y_j} \times V_{x,y_j} \right)
               \subset f^{-1}((-\infty, \tilde r))
 \end{align*}
 for \ $x\in g^{-1}((-\infty,r))$.
\  Especially, given \ $x\in g^{-1}((-\infty,r))$, \  for all \ $x^*\in \bigcap_{i=1}^{k(x)}U_{x,y_i}$ \ and \ $y^*\in Y$ \ we have \ $f(x^*,y^*)<\tilde r$,
 and hence \ $g(x^*) \leq \tilde r <r$ \ for each \ $x^*\in \bigcap_{i=1}^{k(x)}U_{x,y_i}$.
\ From this we can derive
\[
g^{-1}((-\infty,r))=\bigcup_{x\in g^{-1}((-\infty,r))}\left(\bigcap_{i=1}^{k(x)}U_{x,y_i}\right).
\]
On the right hand side of the above equality there is a union of open sets in \ $X$, \ which is open, so $g^{-1}((-\infty,r))$ is open as well.

The  family \ $\{ (-\infty,r), (r,\infty) : r\in\RR\}$ \ constitutes a subbase of the usual topology of $\RR$,
 which implies that the preimage of every open set of $\RR$ by $g$ is open.
Thus $g$ is continuous.
\proofend

\section*{Acknowledgements}
We would like to thank the referees for their comments (especially for the second proof of the first part of Theorem \ref{Pro_random_mean_exp})
 that helped us to improve the paper.

\bibliographystyle{plain}
\bibliography{bb}

\end{document}